\documentclass[11pt,fleqn, a4]{article}

\usepackage[french]{babel}
\usepackage{inputenc} 
\usepackage[T1]{fontenc}

\usepackage{color}
\usepackage{amssymb}
\usepackage{amsmath}

\usepackage{dsfont}

\newcommand{\dcb}{\begin{array}{lll}}
\newcommand{\dce}{\end{array}}
\newcommand{\ebe}{\begin{enumerate}\setlength{\baselineskip}{13pt}\setlength{\parskip}{0pt}}
\newcommand{\dbe}{\end{enumerate}}

\newcommand{\ibegin}{\begin{itemize}\setlength{\baselineskip}{19pt}\setlength{\parskip}{7pt}}
\newcommand{\iend}{\end{itemize}}
\newcommand{\ok}{\rule{4pt}{6pt}}%

\newtheorem{Theorem}{Theorem}[section]
\newtheorem {Cor}{Corollary}[section]
\newtheorem {definition}{Definition}[section]
\newtheorem {pro}{Proposition}[section]
\newtheorem {Lemma}{Lemma}[section]
\newtheorem {rem}{Remark}[section]

\newtheorem {assumption}{Assumption}[section]

\newcommand {\bd}{\begin{definition}}
\newcommand {\ed}{\end{definition}}
\newcommand {\bl}{\begin{Lemma}}
\newcommand {\el}{\end{Lemma}}
\newcommand {\bcor}{\begin{Cor}}
\newcommand {\ecor}{\end{Cor}}
\newcommand {\brem }{\begin{rem} \rm }
\newcommand {\erem }{\end{rem}}
\newcommand{\bethe}{\begin{Theorem}}
\newcommand{\ethe}{\end{Theorem}}
\newcommand {\bassumption}{\begin{assumption}}
\newcommand {\eassumption}{\end{assumption}}

\def\proof{\noindent {\bf Proof. $\, $}}

\def \ind{1\!\!1\!}

\def\cro#1{\langle #1\rangle}

\def\F{{\cal F}}

\def\G{{\cal G}}

\def\ff{{\mathbb F}}

\def\Q{\mathbb Q}

\newcommand{\stocint}{{\centerdot\hspace{1pt}}}
\newcommand{\pv}{{\Gamma\!}}

\newcommand{\ii}{ {\!\raisebox{1pt}{$\dag$}\hspace{1pt}}}
\newcommand{\mrt}{{$\mathfrak{M}$rp}}
\newcommand{\transp}{{^\top\!}}

\date{}
\begin{document}

\title{Martingale representation property in progressively enlarged filtrations}

\author{Monique Jeanblanc
\thanks{This work has benefited financial support by Fédération Bancaire Française, Chaire Risque de Crédit}
\thanks{Monique Jeanblanc. Laboratoire Analyse et Probabilités, Université d'Evry Val D'Essonne, France\newline
\hspace*{17pt}Email:
monique.jeanblanc@univ-evry.fr}\quad  
Shiqi Song\thanks
{Shiqi Song. Laboratoire Analyse et Probabilités, Université d'Evry Val D'Essonne, France\newline 
\hspace*{17pt}Email: shiqi.song@univ-evry.fr} 
}

\maketitle 

\textbf{Abstract}
Consider $\mathbb{G}$ the progressive enlargement of a filtration $\mathbb{F}$ with a random time $\tau$. Assuming that, in $\mathbb{F}$, the martingale representation property holds, we examine conditions under which the martingale representation property holds also in $\mathbb{G}$. A general methodology is developed in this paper, with results covering every known (classical or recent) examples.

\textbf{Kay words} Progressive enlargement of filtration, martingale representation property, honest time, immersion condition, change of probability measures, $(\mathcal{H}')$-hypothesis, credit risk modeling.

\textbf{MSC class.} 60G07, 60G44, 91G40, 97M30.

\

\

\section{Introduction}  
\label{intro}

The theory of progressive enlargement of filtration is a fundamental technique in default risk modeling. In this context, a natural question is whether the progressive enlargement of a filtration possesses the martingale representation property. Many discussions are made on the question (cf. \cite{BJR, BSJ,  CJZ, JC2, Kusuoka, TXY}). On the other hand, as the intensity process is almost the only element that is calibrated from market data, the paper \cite{JS,JS2} studies the problem how to construct models for a given default intensity. It is proved that infinitely many of such models exist and \cite{JS2} provides a procedure (called $\natural$-model in \cite{JS2} and "an evolution model" in this paper) which produces systematically models with the given intensity. (See \cite{song-natural} for a more complete study.) One gets back then to the question if the evolution model satisfies the martingale representation property. It happens that no of the known techniques applies on the evolution model. That motivates the present paper.

We consider a filtration $\mathbb{F}$ in which
the martingale representation property holds: there exists a
multi-dimensional $\mathbb F$-martingale $W$ such that
any $\mathbb{F}$ local martingale null at the origin can be written as a  stochastic integral with respect to $W$. Let $\mathbb{G}$ be the progressive
enlargement of $\mathbb{F}$ with a random time $\tau$. We investigate conditions under which the martingale representation
property also holds in $\mathbb{G}$. 

The question of the martingale representation property in $\mathbb{G}$ had been
considered from the very beginning of the theory of enlargement of
filtration. In the case of an honest time $\tau$, it is proved in
 \cite{barlow}  that if the weak martingale representation property (i.e., a representation with
stochastic integrals with respect to a compensated random measure)
holds in $\mathbb{F}$, it also holds in $\mathbb{G}$. The
structure of $\mathbb{G}$-martingales  is considered  in the same
situation in \cite[Th\'eor\^eme 5.12]{J} where it is proved that
the space of $\mathbb{G}$-square integrable martingales is
generated by two families : the family $\widetilde{\cal {X}}$ of
$\widetilde{X}$, where $\widetilde{X}$ denotes the
$\mathbb{G}$-martingale part of a bounded $\mathbb{F}$-martingale
$X$, and the family $\mathcal {J}$ of the bounded
$\mathbb{G}$-martingales of the form
$v\ind_{[\tau,\infty)}-(v\ind_{[\tau,\infty)})^{\mathbb{G}\cdot
p}$, where $v\in {\cal {G}}_\tau$ and
$(v\ind_{[\tau,\infty)})^{\mathbb{G}\cdot p}$ denotes the
$\mathbb{G}$-predictable dual projection. This study is also
specified in \cite{AJKY} when $\mathbb{F}$ is the natural
filtration of a Brownian motion $B$ and $\tau$ is the end of a
predictable set (hence is honest), avoiding the $\mathbb{F}$
stopping times. The  proof of these results relies on the fact
that the $\mathbb{G}$-predictable sets are generated by the
$\mathbb{F}$-predictable sets and by the random interval
$[0,\tau]$ (cf. \cite{barlow, J}), which is a particular property
of honest time. If we restrict the study on the random interval
$[0,\tau]$, the $\mathbb{G}$-predictable sets coincide with the
$\mathbb{F}$-predictable sets, whatever is the random time $\tau$.
Based on this point, under the technical condition $(C)$ of the
filtration $\mathbb F$, martingales of the form
$K_t=\mathbb{E}[k_\tau|\mathcal{G}_t], t\geq 0$, where $k$ is an
$\mathbb{F}$-predictable process, are studied in \cite{BSJ} and it
is proved that these martingales  $K$ are in the stable space
generated by $\widetilde{\mathcal {X}}$ and $\mathcal{J}$
restricted on $[0,\tau]$ (cf. also \cite{BJR}). Besides the honest
time condition which gives  to the $\mathbb{G}$-predictable sets a
simple structure, the well-known Jacod's criterion (cf. \cite{Jacod}) is another condition
which enables ones to make efficient computations in $\mathbb{G}$.
Actually under Jacod's criterion, $\tau$ behaves much like a
random variable independent of $\mathcal{F}_\infty$. Under this
condition, \cite{JC2} obtain the same generating property of
$\widetilde{\mathcal{X}}$ and $\mathcal{J}$. This idea is further
developed in \cite{CJZ} (and also in \cite{TXY}
for the weak martingale representation property in $\mathbb{G}$). Finally, a very
different and elementary approach is proposed in \cite{Kusuoka}
where the actual martingale representation property is obtained in
$\mathbb{G}$ when $\mathbb{F}$ is a Brownian filtration satisfying
the immersion assumption, i.e., all $\mathbb{F}$-martingales are
$\mathbb{G}$-martingales.

The key point which allows to establish  the martingale
representation property in $\mathbb{G}$ is how the "projection" of
a $\mathbb{G}$-local martingale onto the stable space generated by
$\widetilde{\mathcal {X}}$ can be computed. The contribution of this paper is a general methodology for computing such projections. This methodology gives a unified proof of all the results mentioned before and is applicable to the evolution model of \cite{JS2}. The main consequences of this methodology are presented in Theorem \ref{beforeTau}, Theorem \ref{mrt_after_default}, Theorem \ref{put-together} and Theorem \ref{diez}.

As in the literature, we deal with the problem separately on the
time interval $[0,\tau]$ and on the time interval $(\tau,\infty)$.
For any random time $\tau$, the $\mathbb{G}$-martingales
restricted on $[0,\tau]$  are linked with the
$\mathbb{F}$-martingales in the following formula (see \cite{BJR})
: Let $N^\tau$ be a bounded $\mathbb{G}$-martingale $N$ stopped at
$\tau$. Then, noting $Z_t={\mathbb Q}[t<\tau \vert
\mathcal{F}_{t}]$ and $H_t=\ind_{\{\tau \leq t\}}$
\begin{equation}\label{bdd}
N^\tau_t
=\frac{\mathbb{E}[N_\tau\ind_{\{t<\tau\}}|\mathcal{F}_{t}]}{Z_t}(1-H_t)+N_\tau
H_t,\ 0\leq t <\infty.
\end{equation}
Developing carefully the right hand side of the above identity by
integration by parts formula, we obtain the "projection" of
$N^\tau$ with respect to  $\widetilde{\mathcal {X}}$, with no
supplementary assumption (which generalizes the result in
\cite{BSJ, BJR}). See Section \ref{beforeDefault} for details.

The situation on the time interval $(\tau,\infty)$ is more
complicated. All the  known results are   obtained under extra
conditions. Our approach is based on the method of local solutions
proposed in \cite{song} (cf. \cite{song-local} for a recent account). It has been shown in \cite{song-local}
that the local solution method is  efficient   in the study of the
problem of enlargement of filtration. The well-known Jacod's
criterion and honest time assumptions are two particular cases
where the local solution method is satisfied. Following the
principle of the local solution method, we introduce the notion of
$s\!\mathcal{H}$-measure with covering condition. An
$s\!\mathcal{H}$ measure is a change of probability which creates
locally a relationship between $\mathbb{G}$ and $\mathbb{F}$
similar to the immersion situation considered in \cite{Kusuoka}.
The covering condition says that we can do so at every time point in
$(\tau,\infty)$ (or $(0,\infty)$). This new notion provides a
technique (a mixture of \cite{Kusuoka} and \cite{song}) to compute
the "projection" with respect to $\widetilde{\mathcal {X}}$. See
Section \ref{afterdefault} for details.

In Section \ref{examples}, it is shown that the notion of
$s\!\mathcal{H}$ measure   with covering condition is satisfied in
the situations of the classical results mentioned at the beginning of this section. We provide new proofs of these results. Beyond the classical examples, the notion of
$s\!\mathcal{H}$ measure with covering condition is also
satisfied in situations such as the evolution model in \cite{JS2}, where the
classical assumptions are not available. See Section \ref{JSmodel}. We recall that this paper was initiated by the question if the evolution model in \cite{JS2} satisfies the martingale representation property.

The study of the martingale representation property in
$\mathbb{G}$ provides precise information on the structure of the
family of $\mathbb{G}$-local martingales. It also gives powerful
computing techniques. This study has   applications in
mathematical modeling of financial market. As an example, consider
the notion of market completeness, a notion which guarantees
perfect pricing and hedging. However in general, when the
filtration changes from $\mathbb{F}$ to $\mathbb{G}$, the
completeness property will be lost. We ask questions : What are
the elements which make the loss of the completeness ? How to
quantify the incompleteness ? In what is this incompleteness a bad
thing for the market ? Is it possible to complete an incomplete
market in introducing complementary tradable assets ? As answer,
we   state that it is especially the gap between
$\mathcal{G}_{\tau-}$ and $\mathcal{G}_{\tau}$ which creates   the
  incompleteness. In some cases, this gap can be completed by
introducing  new tradable assets. Theorem
\ref{black-scholes-completion} shows that, if $\mathbb{F}$ is the
filtration of a Black-Scholes model, if $\tau$ is a pure honest
default time, the market with filtration $\mathbb{G}$ is not
complete (it is not possible to hedge the risk associated with the
time $\tau$), but it becomes complete if it is supplemented with
(merely) two defaultable zero-coupon bonds.

\

\section{Preliminary results}
\label{preliminary}

We recall some basic facts in stochastic calculus.

\subsection{$d$-dimensional stochastic integrals}
\label{stocInteg}

A stochastic basis $(\Omega, \mathcal{A},\mathbb{Q},\mathbb{F})$ is a quadruplet, where $(\Omega, \mathcal{A},\mathbb{Q})$ is a probability space and $\mathbb{F}$ is a filtration of sub-$\sigma$-algebras of $\mathcal{A}$, satisfying the usual conditions. 

Given a stochastic basis, we introduce different spaces of martingales. The basic one is $\mathcal{M}({\mathbb{Q}},{\mathbb{F}})$ the space of all $({\mathbb{Q}},{\mathbb{F}})$ martingales. Various subscript or superscript may be used to indicate various derivatives of the space $\mathcal{M}({\mathbb{Q}},{\mathbb{F}})$, especially $\mathcal{M}_{loc,0}({\mathbb{Q}},{\mathbb{F}})$ for local martingales null at the origin, or $\mathcal{M}_{0}^\infty({\mathbb{Q}},{\mathbb{F}})$ for bounded martingales null at the origin. The space $\mathcal{H}^p_0({\mathbb{Q}},{\mathbb{F}})$ is the subspace of $M\in \mathcal{M}_0({\mathbb{Q}},{\mathbb{F}})$ such that $\mathbb{E}[(\sqrt{[M]_\infty})^p]<\infty$ (cf. \cite[Chapter 10]{Yan} and \cite{Jacodlivre}).

Stochastic integral with respect to multi-dimensional local martingale (cf. \cite{cherny, Jacodlivre, JSh}) will be employed in the computations. For a multi-dimensional $({\mathbb{Q}},{\mathbb{F}})$ local martingale $M$, $\mathcal{I}({\mathbb{Q}},{\mathbb{F}},M)$ denotes the family of the multi-dimensional ${\mathbb{F}}$-predictable processes which are $M$-integrable under ${\mathbb{Q}}$ (cf. \cite[Chapitre IV.4 (4.59)]{Jacodlivre} or \cite{cherny}). For $J\in \mathcal{I}({\mathbb{Q}},{\mathbb{F}},M)$, the stochastic integral is denoted by $J\stocint M$ and the space of all such stochastic integrals is denoted by $\mathcal{M}_{loc,0}({\mathbb{Q}},{\mathbb{F}},M)$. Note that, by definition, the stochastic integrals are null at the origin.

We define a topology in $\mathcal{I}({\mathbb{Q}},{\mathbb{F}},M)$. A sequence $(J_n)_{n\geq 1}$ in $\mathcal{I}({\mathbb{Q}},{\mathbb{F}},M)$ is said to converge to an element $J$ in $\mathcal{I}({\mathbb{Q}},{\mathbb{F}},M)$, if there exists a sequence of stopping times $(T_k)_{k\geq 1}$ converging to the infinity such that,  for any $k\geq 1$, the sequence of martingales $(J_n^{T_k}\stocint M)_{n\geq 1}$ converges to $J^{T_k}\stocint M$ in $\mathcal{H}_0^1$. We note that the bounded predictable processes are dense in $\mathcal{I}({\mathbb{Q}},{\mathbb{F}},M)$. 

The notion of stochastic integral is also employed for a multi-dimensional $({\mathbb{Q}},{\mathbb{F}})$ semimartingale $X$ (cf. \cite[Chapter III.6b]{JSh} or \cite{cherny}). Similarly, the space of $X$-integrable processes will be denoted by $\mathcal{I}({\mathbb{Q}},{\mathbb{F}},X)$ and the stochastic integral will be denoted by $J\stocint X$. Recall $(J\stocint X)_0=0$ by definition. Sometimes, the stochastic integral will also be denoted by $(\int_0^t J_s dX_s)_{t\geq 0}$. Recall also that, when $X$ has finite variation, $(J\stocint X)$ is simply the pathwise Lebesgue-Stieltjes integral.

\

\subsection{Stochastic integrals in different filtrations under different probability measures}

A technical problem that we will meet in this paper is that a stochastic integral defined in a filtration under some probability measure will be also considered in another filtration under another probability measure. Assumption \ref{equalintegral} below states conditions which ensure that the variously defined stochastic integrals coincide.

Given a measurable space $(\Omega,\mathcal{A})$, let ${\mathbb{Q}}$ and
${\mathbb{P}}$ be two probability measures on $\mathcal{A}$. Let ${\mathbb{F}}=({\mathcal{F}}_t)_{t\geq
0}$ and ${\mathbb{G}}=({\mathcal{G}}_t)_{t\geq 0}$ be two
right-continuous filtrations in $\mathcal{A}$. 
Let $X$ be a multi-dimensional càdlàg process and  $0\leq S\leq T$ be two given random variables. Consider the following assumption.

\bassumption  \label{a1} 
\ebe 

\item 
$({\mathbb{Q}},{\mathbb{F}})$ (resp.
$({\mathbb{P}},{\mathbb{G}})$) satisfies the usual conditions.

\item 
$S,T$
are ${\mathbb{F}}$-stopping times and
${\mathbb{G}}$-stopping times.

\item
$X$ is
a $({\mathbb{Q}},{\mathbb{F}})$ semimartingale and a
$({\mathbb{P}},{\mathbb{G}})$ semimartingale.

\item   
The probability ${\mathbb{P}}$ is equivalent to
${\mathbb{Q}}$ on
${\mathcal{F}}_\infty\vee{\mathcal{G}}_\infty$.
 
\item   
for any
${\mathbb{F}}$-predictable process $J$, $J\ind_{(S,T]}$ is a
${\mathbb{G}}$-predictable process. 
\dbe 
\eassumption

\bl\label{equalintegral}
Suppose Assumption \ref{a1}. Let $J$ be an ${\mathbb{F}}$-predictable process. Suppose $$
J\ind_{(S,T]}\in\mathcal{I}({\mathbb{Q}},{\mathbb{F}},X)\cap \mathcal{I}({\mathbb{P}},{\mathbb{G}},X).
$$ 
Then, the stochastic integral $J\ind_{(S,T]}\stocint X$ defined in the two senses gives the same process.
\el

\proof
The lemma is true for elementary $\mathbb{F}$ predictable process $J$. Apply \cite[Theorem 1.4]{Yan} (monotone class theorem), \cite[Remark(ii) of Definition 4.8]{cherny} and \cite[Lemma 4.12]{cherny}, we see that the lemma is true for bounded $\mathbb{F}$ predictable process. Again by \cite[Remark(ii) of Definition 4.8]{cherny} and \cite[Lemma 4.11]{cherny}, the lemma is proved. \ok

We can also consider the question of semimartingale decomposition under Assumption \ref{a1}. But we will do it only in an enlarged filtration. Suppose therefore ${\mathcal{F}}_t\subset {\mathcal{G}}_t, t\geq 0$. Suppose that a $({\mathbb{Q}},{\mathbb{F}})$ local martingale $X$ remains a $({\mathbb{Q}},{\mathbb{G}})$ (special) semimartingale. Let $X=M+V$ be the canonical decomposition, where $M$ is a $({\mathbb{Q}},{\mathbb{G}})$ local martingale and $V$ is a ${\mathbb{G}}$ predictable càdlàg process with finite variation. 

\bl\label{jfg}
Let $J\in \mathcal{I}({\mathbb{Q}},{\mathbb{F}},X)$ and denote $Y=J\stocint X$. If $Y$ is a $({\mathbb{Q}},{\mathbb{G}})$ semimartingale, we have $J\in\mathcal{I}({\mathbb{Q}},{\mathbb{G}},M)\cap\mathcal{I}({\mathbb{Q}},{\mathbb{G}},V)$ and $Y=J\stocint M+J\stocint V$.
\el

\proof 
The key point is to prove the $M$-integrabiliy and the $V$-integrability of $J$. Without loss of the generality, suppose $Y\in\mathcal{H}^1_0({\mathbb{Q}},{\mathbb{F}})$. Let $J_n=\ind_{\{\sup_{1\leq i\leq d}|J_i|\leq n\}}J$ for $n\in\mathbb{N}$. By \cite[Corollaire (1.8)]{J}, for some constant $C$, $$
\mathbb{E}[\sqrt{[J_n\stocint M,J_n\stocint M]_\infty}]
\leq
C\mathbb{E}[\sqrt{[J_n\stocint X,J_n\stocint X]_\infty}]
\leq 
C\mathbb{E}[\sqrt{[Y,Y]_\infty}]<\infty, \ n\in\mathbb{N}.
$$
Taking $n\uparrow\infty$, we prove $J\in\mathcal{I}({\mathbb{Q}},{\mathbb{G}},M)$. Because $Y$ is a semimartingale in $\mathbb{G}$, applying Lemma \ref{equalintegral},$$
\dcb
&&\ind_{\{\sup_{1\leq i\leq d}|J_i|\leq n\}}\stocint Y \ \mbox{ (in ${\mathbb{G}}$)}

=\ind_{\{\sup_{1\leq i\leq d}|J_i|\leq n\}}\stocint Y \ \mbox{ (in ${\mathbb{F}}$)}

=J_n\stocint X \ \mbox{ (in ${\mathbb{F}}$)}\\

&=&J_n\stocint X \ \mbox{ (in ${\mathbb{G}}$)}

=
J_n\stocint M +J_n\stocint V \ \mbox{ (in ${\mathbb{G}}$)}.
\dce
$$
As in \cite[Chapter III.6b]{JSh}, we represent the components $V_i$ of $V$ in the form $a_i\stocint F$, where the processes $a_i$ and $F$ are supposed to be ${\mathbb{G}}$ predictable and $F$ is càdlàg increasing. We have$$
\dcb
&&|\sum_{i=1}^dJ_ia_i|\ind_{\{\sup_{1\leq i\leq d}|J_i|\leq n\}}\stocint F\\
&=&
\mbox{sgn}(\sum_{i=1}^dJ_ia_i)\ind_{\{\sup_{1\leq i\leq d}|J_i|\leq n\}}J\stocint V\\
&=&\mbox{sgn}(\sum_{i=1}^dJ_ia_i)\ind_{\{\sup_{1\leq i\leq d}|J_i|\leq n\}}\stocint Y
-
\mbox{sgn}(\sum_{i=1}^dJ_ia_i)\ind_{\{\sup_{1\leq i\leq d}|J_i|\leq n\}}J\stocint M \ \mbox{ (in ${\mathbb{G}}$)}.
\dce
$$
As $Y$ and $J\stocint M$ are semimartingales in $\mathbb{G}$, by \cite[Remark(ii) of Definition 4.8]{cherny} and \cite[Lemma 4.11 and Lemma 4.12]{cherny}, the terms on the right hand side of this equality converge in probability when $n\uparrow \infty$. Consequently, $J\in\mathcal{I}({\mathbb{Q}},{\mathbb{G}},V)$ (cf. \cite[Definition 3.7]{cherny}) and the lemma follows. \ok

\

\subsection{Martingale representation property}
\label{recallmrt}

In this subsection, we give the definition of the martingale representation property that we adopt in this paper, and we recall some related results with short proofs.

Given a stochastic basis $(\Omega, \mathcal{A},\mathbb{Q},\mathbb{F})$, let $W$ be a multi-dimensional ${\mathbb{F}}$-adapted càdlàg process. We say that $W$ has the martingale representation property in the filtration ${\mathbb{F}}$ under the probability ${\mathbb{Q}}$, if $W$ is a $({\mathbb{Q}},{\mathbb{F}})$ local martingale, and if $$
\mathcal{M}_{loc,0}({\mathbb{Q}},{\mathbb{F}},W) = \mathcal{M}_{loc,0}({\mathbb{Q}},{\mathbb{F}}).
$$
We will call the process $W$ the driving process. The martingale representation property will be denoted by $\mathfrak{M}$rp(${\mathbb{Q}},{\mathbb{F}},W$), or simply by \mrt.

We introduce the operator\footnote{Note that the operator $\ii$
depends on $({\mathbb{Q}},{\mathbb{F}})$. The context in
which $\ii$ is used will help to avoid the ambiguity.}
 $\ii$ : For $\zeta$ ${\mathcal{F}}_\infty$-measurable ${\mathbb{Q}}$-integrable random variable, we denote by $\zeta^\dag$ the martingale$$
\zeta^\dag_t=	{\mathbb{Q}}[\zeta|{\mathcal{F}}_t]-{\mathbb{Q}}[\zeta|{\mathcal{F}}_0], t\geq 0.
$$

\bl\label{boundedenough}
Suppose that $W$ is a multi-dimensional $({\mathbb{Q}},{\mathbb{F}})$ local martingale. Let $\mathcal{C}$ be a $\pi$-class such that $\sigma(\mathcal{C})={\mathcal{F}}_\infty$. If, for any $A\in\mathcal{C}$, the 
$({\mathbb{Q}},{\mathbb{F}})$-martingale\hspace{3pt} $(\ind_A)^\dag$ is an element in $\mathcal{M}_{loc,0}({\mathbb{Q}},{\mathbb{F}},W)$, then \mrt$({\mathbb{Q}},{\mathbb{F}},W)$ holds. 

\el

\proof We know that, if $\mathcal{M}_0^\infty({\mathbb{Q}},{\mathbb{F}})\subset\mathcal{M}_{loc,0}({\mathbb{Q}},{\mathbb{F}},W)$, then the space $\mathcal{H}^1_0({\mathbb{Q}},{\mathbb{F}})$ is contained in $\mathcal{M}_{0}({\mathbb{Q}},{\mathbb{F}},W)$  (cf. \cite[Theorem 10.5]{Yan}) and consequently \mrt$({\mathbb{Q}},{\mathbb{F}},W)$ holds. Let $\Pi$ be the family of bounded ${\mathcal{F}}_\infty$-measurable random variables $\zeta$ such that the $({\mathbb{Q}},{\mathbb{F}})$-martingales $\zeta^\dag$ belong to $\mathcal{M}_{loc,0}({\mathbb{Q}},{\mathbb{F}},W)$. Applying the monotone class theorem, we see that the space $\Pi$ contains all bounded $\sigma(\mathcal{C})={\mathcal{F}}_\infty$ measurable random variables, which means exactly $\mathcal{M}_0^\infty({\mathbb{Q}},{\mathbb{F}})\subset\mathcal{M}_{loc,0}({\mathbb{Q}},{\mathbb{F}},W)$. \ok

Other characteristic conditions for \mrt\ exist.

\bl\label{caracMRT}
Let $W$ be a $d$-dimensional $({\mathbb{Q}},{\mathbb{F}})$ local martingale with $W_0=0$. The following statements are equivalent :
\ebe
\item
$\mathcal{M}_{loc,0}({\mathbb{Q}},{\mathbb{F}},W)=\mathcal{M}_{loc,0}({\mathbb{Q}},{\mathbb{F}})$, i.e., \mrt(${\mathbb{Q}},{\mathbb{F}},W$) holds.
\item
Any $L\in \mathcal{M}_{loc,0}({\mathbb{Q}},{\mathbb{F}})$ such that $LW_i\in \mathcal{M}_{loc,0}({\mathbb{Q}},{\mathbb{F}})$ for all $1\leq i\leq d$,  is null.
\item
Any $L\in \mathcal{M}^\infty_{0}({\mathbb{Q}},{\mathbb{F}})$ such that $LW_i\in \mathcal{M}_{loc,0}({\mathbb{Q}},{\mathbb{F}})$ for all $1\leq i\leq d$,  is null.
\dbe 
\el

\proof This is the consequence of \cite[Corollaire(4.12) and Proposition(4.67)]{Jacodlivre}\ok

We recall below that the property \mrt\ is invariant by change of probabilities. (However, the driving process may change.) Suppose \mrt(${\mathbb{Q}},{\mathbb{F}},W$). For a probability measure ${\mathbb{P}}$ locally equivalent to ${\mathbb{Q}}$, set$$
\eta_t=\left.\frac{d{\mathbb{P}}}{d{\mathbb{Q}}}\right|_{{\mathcal{F}}_t}, \ 0\leq t<\infty.
$$
The process $\eta$ is a strictly positive $({\mathbb{Q}},{\mathbb{F}})$-martingale. We have the following result.

\bl\label{Girsanovrepresentation} 
Suppose that the predictable brackets $\cro{\eta,W_i}, 1\leq i\leq d,$ exist in ${\mathbb{F}}$ under ${\mathbb{Q}}$ and \mrt(${\mathbb{Q}},{\mathbb{F}},W$) holds. Set 
\begin{equation}\label{weta}
W^{[\eta]}_{i,t}=W_{i,t}-\int_0^t\frac{1}{\eta_{s-}} d\cro{\eta,W_i}_s, 0\leq t<\infty, 1\leq i\leq d.
\end{equation}
Then, we have \mrt(${\mathbb{P}},{\mathbb{F}},W^{[\eta]}$). 
\el

\proof Note that, by Girsanov's theorem, the process $W^{[\eta]}$ is a $({\mathbb{P}},{\mathbb{F}})$ local martingale. By Lemma \ref{caracMRT}, it is enough to prove that, for $L\in\mathcal{M}_{loc,0}({\mathbb{P}},{\mathbb{F}})$, if $LW^{[\eta]}_i\in\mathcal{M}_{loc,0}({\mathbb{P}},{\mathbb{F}})$ for $1\leq i\leq d$, the process $L$ is null. But then, if $Y=L\eta$, the two processes $Y$ and $YW^{[\eta]}_i$ belong to $\mathcal{M}_{loc,0}({\mathbb{Q}},{\mathbb{F}})$. This yields that $\cro{Y,W_i}$ exists (under ${\mathbb{Q}}$) and $
\cro{Y-\frac{Y_-}{\eta_-}{\stocint}\eta,W_i}=0, 1\leq i\leq d.
$
Applying Lemma \ref{caracMRT} with the \mrt(${\mathbb{Q}},{\mathbb{F}},W$) property, we conclude $Y=0$, and consequently $L=0$. \ok

\

\subsection{Progressive enlargement of filtrations and $(\mathcal{H}')$ hypothesis}\label{framework}

Given a stochastic basis $(\Omega,\mathcal{A},\mathbb{F},\mathbb{Q})$, let $\tau$ be an $\mathcal{A}$-measurable random variable taking values in $[0,\infty]$ and $\mathbb{G}$ be the progressive enlargement of the filtration $\mathbb{F}$ by the random time $\tau$, i.e., $\mathbb{G}=(\mathcal{G}_t)_{t\geq 0}$ with $\G_t=\cap_{s>t}(\F_{s}\vee \sigma (\tau \wedge s))$ augmented by $(\mathbb{Q},\mathcal{F}_\infty\vee\sigma(\tau))$ negligible sets.  

The $(\mathbb{Q},\mathbb{F})$ optional projection of the process $\ind_{[0,\tau)}$, denoted by $Z$, is a bounded and non negative $(\mathbb{Q},\mathbb{F})$ supermartingale. We denote by $Z=M-A$ its Doob-Meyer's decomposition, where $M$ is a $(\mathbb{Q},\mathbb{F})$ martingale and $A$ is a càdlàg $(\mathbb{Q},\mathbb{F})$-predictable increasing process with $A_0=0$. We note that $A$ coincides with the $(\mathbb{Q},\mathbb{F})$ predictable dual projection of the process $\ind_{\{\tau>0\}}\ind_{[\tau,\infty)}$. We consider also the $(\mathbb{Q},\mathbb{F})$ optional dual projection of the process $\ind_{\{\tau>0\}}\ind_{[\tau,\infty)}$, that we denote by $\hat{A}$. Let $H:=\ind_{[\tau,\infty)}$. According to \cite[Remarques(4.5) 3)]{J}, the process \begin{equation}\label{L-martingale}
L_t=\ind_{\{\tau>0\}}H_t - \int_0^{t\wedge \tau}\frac{dA_s}{Z_{s-}}, \ t\geq 0,
\end{equation} 
is a $(\mathbb{Q},\mathbb{G})$ local martingale. We introduce the following assumption.  

\bassumption\label{HH}
\noindent\textbf{$(\mathcal{H}')$ Hypothesis.} There exists a map $\Gamma$ from $\mathcal{M}_{loc}(\mathbb{Q},\mathbb{F})$ into the space of càdlàg $\mathbb{G}$-predictable processes with finite variation, such that, for any $X\in \mathcal{M}_{loc}(\mathbb{Q},\mathbb{F})$, $\pv(X)_0=0$ and $\widetilde{X}:=X-\pv(X)$ is a $(\mathbb{Q},\mathbb{G})$ local martingale. The operator $\pv$ will be called the drift operator.
\eassumption

Note that, when a $(\mathbb{Q},\mathbb{F})$ local martingale $X$ is a $(\mathbb{Q},\mathbb{G})$-semimartingale, it is a special semimartingale. Consequently, the drift operator $\pv(X)$ is well defined. Many facts are known about the drift operator $\pv$ (see \cite{J}).
In particular, $\ind_{(0,\tau]}{\stocint} \pv(X)$ takes the form :
\begin{equation}\label{formula_before_default}
\ind_{(0,\tau]}{\stocint} \pv(X) = \ind_{(0,\tau]}\frac{1}{Z_{-}}{\stocint} \cro{M+\hat{A}-A,X},
\end{equation} 
where the bracket $\cro{\cdot}$ is computed in the filtration $\mathbb{F}$. The computation of the drift operator $\pv(X)$ on the time interval $(\tau,\infty)$ is more complicated. There do not exist general results. For a long time, the only known example was the case of a honest time $\tau$ for which the drift operator on $(\tau,\infty)$ is given by:
\begin{equation}\label{formula_after_default}
\ind_{(\tau,\infty)}{\stocint} \pv(X) = \ind_{(\tau,\infty)}\frac{1}{1-Z_{-}}{\stocint} \cro{M+\hat{A}-A,X}
\end{equation}
(see \cite{barlow}, \cite{JY2}). Recent researches (cf. \cite{CJZ,ejj,JC,JS,JS2}) show that, unlike the situation on the interval $[0,\tau]$, the operator $\ind_{(\tau,\infty)}{\stocint} \pv(X)$ can have various different forms. 

\

\subsection{Miscellaneous notes}\label{notations}

When we say that a filtration $\mathbb{F}=(\mathcal{F}_t)_{t\geq 0}$ satisfies the usual condition under a probability $\mathbb{P}$, we consider only the $(\mathbb{P},\mathcal{F}_\infty)$-negligible sets in the usual condition.

Relations between random variables are to be understood to be almost sure relations. For a random variable $X$ and a $\sigma$-algebra $\mathcal{F}$, the expression $X\in\mathcal{F}$ means that $X$ is $\mathcal{F}$-measurable. 

We will compute expectations with respect to different probability measures $\mathbb{P}$. To simplify the notation, we will denote the expectation by $\mathbb{P}[\cdot]$ instead of $\mathbb{E}_{\mathbb{P}}[\cdot]$, and the conditional expectation by $\mathbb{P}[\cdot|\mathcal{F}_t]$ instead of $\mathbb{E}_{\mathbb{P}}[\cdot|\mathcal{F}_t]$. 

For any non negative random variable $\rho$, $\mathcal{F}_\rho$ (resp. $\mathcal{F}_{\rho-}$) denotes the $\sigma$-algebra generated by the random variables $U_\rho\ind_{\{\rho<\infty\}}+\xi\ind_{\{\rho=\infty\}}$ where $U$ runs over the family of the $\mathbb{F}$-optional processes (respectively $\mathbb{F}$-predictable processes) and $\xi\in\mathcal{F}_\infty$

Let $D$ be a subset of $\Omega$ and $\mathcal{T}$ be a $\sigma$-algebra on $\Omega$. We denote by $D\cap\mathcal{T}$ the family of all subsets $D\cap A$ with $A$ running through $\mathcal{T}$. If $D$ itself is an element in $\mathcal{T}$, $D\cap\mathcal{T}$ coincides with $\{A\in\mathcal{T}: A\subset D\}$. Although $D\cap\mathcal{T}$ is a $\sigma$-algebra on $D$, we will use it mainly as a family of subsets of $\Omega$. We use the symbol "+" to present the union of two disjoint subsets. For two disjoint sets $D_1,D_2$ in $\Omega$, and two families $\mathcal{T}_1,\mathcal{T}_2$ of sets in $\Omega$, we denote by $D_1\cap\mathcal{T}_1+D_2\cap\mathcal{T}_2$ the family of sets $D_1\cap B_1+D_2\cap B_2$ where $B_1\in\mathcal{T}_1,B_2\in\mathcal{T}_2$. For a probability $\mathbb{P}$, we say $D\cap \mathcal{T}_1=D\cap\mathcal{T}_2$ under $\mathbb{P}$, if, for any $A\in\mathcal{T}_1$, there exists a $B\in\mathcal{T}_2$ such that $(D\cap A)\Delta(D\cap B)$ is $\mathbb{P}$-negligible, and vice versa. 

The jump at the origin of a càdlàg process $X$ is by definition : $\Delta_0 X=X_0$.

\

\section{\mrt\ property before $\tau$}
\label{beforeDefault}

In this section, we assume the setting of the subsection \ref{framework} with the process $Z,A$ and the notation $\widetilde{X}:=X-\pv(X)$. We need the following lemma from \cite{yorlemma}

\bl\label{Z-positive}
We define $Z_{0-}=1$. Let $R=\inf\{t\geq 0: Z_{t-}=0 \mbox{ or } Z_t=0\}$ and, for $n\geq 1$, $R_n=\inf\{t\geq 0: Z_{t-}<\frac{1}{n} \mbox{ or } Z_t<\frac{1}{n}\}$. Then, $\sup_{n\geq 1}R_n=R$, $\tau\leq R$ and $Z_{\tau-}>0$ on $\{\tau<\infty\}$. 
\el

We introduce

\bassumption
\label{trois-hypos}{ }\
\ebe
\item[(i)]
\mrt$(\mathbb{Q},\mathbb{F},W)$ holds for a $d$-dimensional $\mathbb{F}$-adapted càdlàg process $W$.
\item[(ii)]
$(\mathcal{H}')$-hypothesis holds between $\mathbb{F}$ and $\mathbb{G}$.  
\dbe
\eassumption

We use the notations $\widetilde{X}=X-\pv(X)$ (introduced in Assumption \ref{HH}) and $H=\ind_{[\tau,\infty)}$. Here is the main result in this section, which generalizes \cite[Lemme(5.15)]{J}, \cite[Theorem 1]{BSJ} and \cite[Theorem 3.3.2]{BJR} :

\bethe\label{beforeTau}
Suppose Assumption \ref{trois-hypos}. Then, for any bounded $\zeta\in\mathcal{G}_\tau$, there exist an $\mathbb{F}$-predictable process $J$ such that $J\ind_{[0,\tau]}\in\mathcal{I}(\mathbb{Q},\mathbb{G},\widetilde{W})$ and 
\begin{equation}\label{mrp_before_defaut}
\mathbb{Q}[\zeta|\mathcal{G}_t]=\mathbb{Q}[\zeta|\mathcal{G}_0]+\int_0^t(1-H_{s-})J_sd\widetilde{W}_s\\
+\ind_{\{\tau>0\}}(\zeta-X_{\tau}) H_t-\int_0^t K_s(1-H_{s-})\frac{1}{Z_{s-}}dA_s,
\end{equation}
for $t\geq 0$, where $K$ is a bounded $\mathbb{F}$-predictable process such that $K_\tau\ind_{\{0<\tau<\infty\}}=\mathbb{Q}[(\zeta-X_\tau)\ind_{\{0<\tau<\infty\}}|\mathcal{F}_{\tau-}]$, and the process $X$ is defined as$$
X_t=\frac{\mathbb{Q}[\zeta\ind_{\{t<\tau\}}|\mathcal{F}_{t}]}{Z_t}\ind_{\{t<R\}},\ 0\leq t<\infty.
$$

\ethe

\proof The proof could be build from \cite[Lemma (5.15)]{J}. But we prefer a proof based on a direct computation of the martingale $\mathbb{Q}[\zeta|\mathcal{G}_t], t\geq 0$.  

Consider a bounded random variable $\zeta\in \mathcal{G}_{\tau}$. From \cite{BJR, JR}, we have the decomposition formula$:$ $$
\mathbb{Q}[\zeta|\mathcal{G}_{t}]\\
=\frac{\mathbb{Q}[\zeta\ind_{\{t<\tau\}}|\mathcal{F}_{t}]}{Z_t}(1-H_t)+\zeta H_t,\ 0\leq t <\infty.
$$
Since $\tau\leq R$ and $Z_{\tau-}>0$, the process $X$ is well defined, bounded and c\`adl\`ag. By \cite[Lemma 5.2]{barlow}, for any $n\geq 1$, $X^{R_n}$ is a $(\mathbb{Q},\mathbb{F})$ semimartingale, i.e., the process $X$ is a bounded $(\mathbb{Q},\mathbb{F})$ semimartingale on the $\mathbb{F}$-predictable set $\mathtt{B}=\cup_{n\geq 1}[0,R_n]$ (see \cite[Definition 8.19]{Yan}). Note that $\sup_{n\geq 1}R_n=R$ which implies $\cup_{n\geq 1}[0,R_n]\supset [0,R)$. Denote the $(\mathbb{Q},\mathbb{F})$ canonical decomposition on $\mathtt{B}$ of $X$ by $X^{[m]}+X^{[v]}$, where $X^{[m]}$ is a $(\mathbb{Q},\mathbb{F})$ local martingale on $\mathtt{B}$ and $X^{[v]}$ is an $\mathbb{F}$-predictable process with finite variation on $\mathtt{B}$. By the property \mrt($\mathbb{Q},\mathbb{F},W$), there exists a $d$-dimensional $\mathbb{F}$-predictable process $J$ defined on $\mathtt{B}$ such that, for any $n\geq 1$, $J\ind_{[0,R_n]}\in\mathcal{I}(\mathbb{Q},\mathbb{F},W)$ and 
$$
X^{[m]}_t=X^{[m]}_0+\int_0^t J_sdW_s,\ 0\leq t<R.
$$
Since, for every $n\geq 1$, ${J\ind_{[0,R_n]}}\centerdot W$ and $W$ are $(\mathbb{Q},\mathbb{G})$ semimartingales by $(\mathcal{H}')$-Hypothesis, according to Lemma \ref{jfg}, the process $J\ind_{[0,R_n]}$ belongs to $\mathcal{I}(\mathbb{Q},\mathbb{G},\widetilde{W})$ and to $\mathcal{I}(\mathbb{Q},\mathbb{G},\pv(W))$. Lemma \ref{equalintegral} is applicable to $J\ind_{[0,R_n]}$ between $\mathbb{F}$ and $\mathbb{G}$.

Note $\mathcal{G}_{\tau-}=\mathcal{F}_{\tau-}$ (cf. \cite[Lemme (4.4)]{J}). Let $K$ be a bounded $\mathbb{F}$-predictable process such that $$
K_\tau\ind_{\{0<\tau<\infty\}}=\mathbb{Q}[(\zeta-X_\tau)\ind_{\{0<\tau<\infty\}}|\mathcal{G}_{\tau-}]
=\mathbb{Q}[(\zeta-X_\tau)\ind_{\{0<\tau<\infty\}}|\mathcal{F}_{\tau-}].
$$ 
We compute now the martingale $\mathbb{Q}[\zeta|\mathcal{G}_{t}]$ $:$ for $n\geq 1$, for $0\leq t\leq R_n$ 
$$
\dcb
\mathbb{Q}[\zeta|\mathcal{G}_{t}]
&=&X_t(1-H_t)+\zeta H_t\\
&=&\int_0^t(1-H_{s-})dX_s-\int_0^tX_{s-} dH_s + [X, 1-H]_t +\zeta H_t\\

&=&X_0(1-H_0)+\int_0^t(1-H_{s-})J_sdW_s+\int_0^t(1-H_{s-})dX^{[v]}_s+(\zeta-\ind_{\{\tau>0\}}X_{\tau}) H_t\\

&=&\mathbb{Q}[\zeta|\mathcal{G}_{0}]+\int_0^t(1-H_{s-})J_sd\widetilde{W}_s+\int_0^t(1-H_{s-})J_sd\pv(W)_s+\int_0^t(1-H_{s-})dX^{[v]}_s\\
&&+\ind_{\{\tau>0\}}(\zeta-X_{\tau}) H_t-\int_0^t K_s\ind_{\{s\leq\tau\}}\frac{1}{Z_{s-}}dA_s
+\int_0^t K_s\ind_{\{s\leq\tau\}}\frac{1}{Z_{s-}}dA_s.\\
\dce
$$
We note that, since $\int_0^t K_s\ind_{\{s\leq\tau\}}\frac{1}{Z_{s-}}dA_s$ is the $(\mathbb{Q},\mathbb{G})$ predictable dual projection of $\ind_{\{\tau>0\}}(\zeta-X_{\tau}) H$ (see \cite{J}), the process$$
\ind_{\{\tau>0\}}(\zeta-X_{\tau}) H_t-\int_0^t K_s\ind_{\{s\leq\tau\}}\frac{1}{Z_{s-}}dA_s
$$
is a $(\mathbb{Q},\mathbb{G})$ local martingale. This implies that the following $\mathbb{G}$-predictable process with finite variation $$
\int_0^t(1-H_{s-})J_sd\pv(W)_s+\int_0^t(1-H_{s-})dX^{[v]}_s+\int_0^t K_s\ind_{\{s\leq\tau\}}\frac{1}{Z_{s-}}dA_s
$$
is a $(\mathbb{Q},\mathbb{G})$ local martingale, so that it is null. Consequently,
\begin{equation}\label{HJK}
\mathbb{Q}[\zeta|\mathcal{G}_{t}]
=\mathbb{Q}[\zeta|\mathcal{G}_0]+\int_0^t(1-H_{s-})J_sd\widetilde{W}_s +\ind_{\{\tau>0\}}(\zeta-X_{\tau}) H_t-\int_0^t K_s\ind_{\{s\leq\tau\}}\frac{1}{Z_{s-}}dA_s
\end{equation}
for $0\leq t\leq R_n$, $n\geq 1$. 

Let us prove that $\mathbb{Q}[\forall n\geq 1, R_n<\tau]=0$. Actually, on the set $\{\forall n\geq 1, R_n<\tau\}$, there can exist two situations. Firstly the sequence $(R_n)$ is stationary, i.e., for some $n$, $R=R_n<\tau$. It is impossible because $Z_{R+\epsilon}=0,\forall \epsilon>0$ (see \cite[Theorem 2.62]{Yan}) whilst $Z_{\tau-}>0$. When the sequence $(R_n)$ is not stationary, we must have $R=\lim_nR_n=\tau$ and $Z_{R-}=0$. Once again it is impossible because $Z_{\tau-}>0$. 

We conclude that $[0,\tau]\subset \mathtt{B}$. Hence, the process $J$ is well defined on $[0,\tau]$ and the above formula (\ref{HJK}) is true on $[0,\tau]$.
Note that $\mathbb{Q}[\zeta|\mathcal{G}_{t}]$ and $\ind_{\{\tau>0\}}(\zeta-X_{\tau}) H_t-\int_0^t K_s\ind_{\{s\leq\tau\}}\frac{1}{Z_{s-}}dA_s$ are $(\mathbb{Q},\mathbb{G})$ local martingales on the whole $\mathbb{R}_+$. Consequently, $J\ind_{[0,R_n]}$ is in $\mathcal{I}(\mathbb{Q},\mathbb{G},\widetilde{W})$ uniformly for all $n\geq 1$, which implies $J\ind_{[0,\tau]}\in\mathcal{I}(\mathbb{Q},\mathbb{G},\widetilde{W})$. The theorem is proved. \ok

\brem\label{Ldef}
Let $L = \ind_{\{\tau>0\}}\ind_{[\tau,\infty)} - \ind_{(0,\tau]}\frac{1}{Z_{-}}{\stocint} A$. Then, the formula (\ref{mrp_before_defaut}) can be written as
\begin{equation}
\dcb
\mathbb{Q}[\zeta|\mathcal{G}_t]	
&=&\mathbb{Q}[\zeta|\mathcal{G}_0]+\int_0^t\ind_{\{0<s\leq\tau\}}J_sd\widetilde{W}_s
\\
&&\hspace{1.3cm}+\int_0^tK_s\ind_{\{0<s\leq\tau\}}dL_s+\ind_{\{\tau>0\}}(\zeta-X_{\tau}-K_\tau) H_t.
\dce
\end{equation}
This shows that, for any bounded $(\mathbb{Q},\mathbb{G})$-martingale $Y$, the stopped martingale $Y^\tau$ is the sum of a stochastic integral with respect to $(\widetilde{W},L)$ and a $(\mathbb{Q},\mathbb{G})$-martingale of the form $\xi H$, where $\xi\in\mathcal{G}_\tau$ such that $\mathbb{Q}[\xi|\mathcal{G}_{\tau-}]=0$. This remark links the formula (\ref{mrp_before_defaut}) to \cite[Théorème(5.12)]{J}.
\erem

Suppose now that $\{0<\tau<\infty\}\cap\mathcal{G}_\tau=\{0<\tau<\infty\}\cap\mathcal{G}_{\tau-}$ under $\mathbb{Q}$. In this case, $\ind_{\{0<\tau<\infty\}}K_\tau= \ind_{\{0<\tau<\infty\}}(\zeta-X_{\tau})$ and formula (\ref{mrp_before_defaut}) writes as
\begin{equation}\label{MWL}
\dcb
\mathbb{Q}[\zeta|\mathcal{G}_{t}]
&=&\mathbb{Q}[\zeta|\mathcal{G}_0]+\int_0^t\ind_{\{0<s\leq\tau\}}J_sd\widetilde{W}_s +\int_0^t K_s\ind_{\{0<s\leq\tau\}}dL_s.
\dce
\end{equation}
From this identity, applying Lemma \ref{boundedenough}, we conclude that the property \mrt$(\mathbb{Q},\mathbb{G}^{\tau},(\widetilde{W}^\tau,L))$ holds, where $\mathbb{G}^\tau$ denotes the stopped filtration $(\mathcal{G}_{t\wedge \tau}: t\geq 0)$, and $\widetilde{W}^{\tau}$ is $\widetilde{W}$ stopped at $\tau$. We can state this conclusion in another way :

\bethe\label{companion}
Suppose Assumption \ref{trois-hypos}. Then, \mrt$(\mathbb{Q},\mathbb{G}^{\tau},(\widetilde{W}^\tau,L))$  and $\ind_{\{0<\tau<\infty\}}W_\tau\in\mathcal{G}_{\tau-}$ hold if and only if $\{0<\tau<\infty\}\cap\mathcal{G}_\tau=\{0<\tau<\infty\}\cap\mathcal{G}_{\tau-}$ under $\mathbb{Q}$. 
\ethe

\proof We have explained that the condition is sufficient. Suppose now \mrt$(\mathbb{Q},\mathbb{G}^{\tau},(\widetilde{W}^{\tau},L))$ and $\ind_{\{0<\tau<\infty\}}W_\tau\in\mathcal{G}_{\tau-}$. Any $(\mathbb{Q},\mathbb{G}^\tau)$ bounded martingale $Y$ with $Y_0=0$ has a representation as in (\ref{MWL}). We have, therefore, $$
\dcb
&&\ind_{\{0<\tau<\infty\}}\Delta_\tau Y \\
&=& J_\tau\ind_{\{0<\tau<\infty\}}\Delta_\tau\widetilde{W} + \ind_{\{0<\tau<\infty\}}K_\tau(1-\frac{1}{Z_{\tau-}}\Delta_\tau A) \mbox{ (computing $\Delta_\tau Y$ with \cite[Proposition(4.68)]{Jacodlivre})}\\
&=& J_\tau\ind_{\{0<\tau<\infty\}}(W_\tau - W_{\tau-}) - J_\tau\ind_{\{0<\tau<\infty\}}\Delta_\tau\pv(W) + \ind_{\{0<\tau<\infty\}}K_\tau(1-\frac{1}{Z_{\tau-}}\Delta_\tau A),\\

\dce
$$
where $J,K,Z_-,\Delta\pv(W),\Delta A$ are $\mathbb{F}$-predictable processes and $\{0<\tau<\infty\}\in\mathcal{G}_{\tau-}$. Consequently, $\ind_{\{0<\tau<\infty\}}Y_\tau = \ind_{\{0<\tau<\infty\}}Y_{\tau-}+\ind_{\{0<\tau<\infty\}}\Delta_\tau Y \in\mathcal{G}_{\tau-}$, and (cf. \cite[Chapter XX section 22]{DM})$$
\dcb
&&\{0<\tau<\infty\}\cap\mathcal{G}_\tau \\
&=& \{0<\tau<\infty\}\cap\sigma\{\tau, \ind_{\{0<\tau<\infty\}}Y_\tau : \mbox{ $Y$ is a $(\mathbb{Q},\mathbb{G}^\tau)$ bounded martingale}\}\\

&\subset& \{0<\tau<\infty\}\cap\mathcal{G}_{\tau-}.
\dce
$$
The inverse inclusion $\mathcal{G}_{\tau-}\subset \mathcal{G}_{\tau}$ is obvious. The theorem is proved. \ok

\brem
The condition $\ind_{\{0<\tau<\infty\}}W_{\tau}\in\mathcal{G}_{\tau-}$ is satisfied if $W$ is continuous. It is also satisfied if the random time $\tau$ avoids the $\mathbb{F}$ stopping times, since then $W_\tau=W_{\tau-}\in\mathcal{G}_{\tau-}$. 
\erem

\

\section{\mrt\ property after the default time $\tau$}
\label{afterdefault}

We will work on the \mrt\ property on the time horizon $(\tau,\infty)$ according to the local solution methodology developed in \cite{song,song-local}. Roughly speaking, this amounts to proceed the search of \mrt\ property in two steps. Firstly for any point $a\in(\tau,\infty)$, we look for, if possible, a neighbourhood $(s_a,t_a)$ of $a$ and a probability measure $\mathbb{Q}_a$ equivalent to $\mathbb{Q}$ such that under $\mathbb{Q}_a$, the "restriction of the filtration $\mathbb{G}$" on $(s_a,t_a)$ has a "local immersion property on $(s_a,t_a)$", which will ensure an "\mrt\ property on $(s_a,t_a)$" (a local solution). Then, we establish an \mrt\ property on the union $\cup_{a\in(\tau,\infty)}(s_a,t_a)$ (a global solution). This will, in particular, give a solution to our problem when $\cup_{a\in(\tau,\infty)}(s_a,t_a)\supset[\tau,\infty)$. In this section, the above mentioned concepts such as "local immersion property" or "restriction of the filtration $\mathbb{G}$" will be firmly defined and the two steps program will be realized.

\subsection{Fragments of the filtration $\mathbb{G}$}
\label{frags}  

We work in the setting of the subsection \ref{framework}. The following definition represents the notion of the "restriction of the filtration $\mathbb{G}$" on a random interval.

\bd \label{filtrationGST}
For any $\mathbb{G}$-stopping time $T$, we define the $\sigma$-algebra$$
\mathcal{G}^\ast_{T}=\{T<\tau\}\cap\mathcal{F}_{T}+\{\tau\leq T<\infty\}\cap(\sigma(\tau)\vee\mathcal{F}_{T})+\{T=\infty\}\cap(\sigma(\tau)\vee\mathcal{F}_{\infty}).
$$For two $\mathbb{G}$-stopping times $S,T$ such that $S\leq T$, we define the family $\mathbb{G}^{(S,T]}$ of $\sigma$-algebras : $$
\mathcal{G}^{(S,T]}_t 
= \{\{T\leq S\vee t\}\cap A + \{S\vee t<T\}\cap B : A\in\mathcal{G}^\ast_T,\ B\in\mathcal{G}_{S\vee t}\},\ 0\leq t<\infty.
$$
For general $\mathbb{G}$-stopping times $S,T$, we define the family $\mathbb{G}^{(S,T]}= \mathbb{G}^{(S,S\vee T]}$. 

\label{SvT}
For any process $X$, for $\mathbb{G}$-stopping times $S,T$ such that $S\leq T$, we denote $X^{(S,T]}_t=X_{(S\vee t)\wedge T}-X_S, 0\leq t<\infty$. For general $\mathbb{G}$-stopping times $S,T$, we define $X^{(S,T]}=X^{(S,S\vee T]}$. 
\ed

The following results exhibit the properties of the family $\mathbb{G}^{(S,T]}$ in relation with the filtration $\mathbb{G}$. The proofs will be given in the Appendix. 

\bl\label{g-star}
For any $\mathbb{G}$-stopping times $S,T$, we have $\mathcal{G}_{T-}\subset \mathcal{G}^\ast_{T}\subset \mathcal{G}_{T}$ and $$
\mathcal{G}^*_{S\vee T}=\{S<T\}\cap \mathcal{G}^*_{T}+\{T\leq S\}\cap \mathcal{G}^*_{S}.
$$
\el

\begin{pro}  \label{appendix} 
We consider two $\mathbb{G}$-stopping times $S,T$ such that $S\leq T$. We have
\ebe
\item[(1)]  
$\mathbb{G}^{(S,T]}$ is a right-continuous filtration.
\item[(2)]
For a $\mathbb{G}^{(S,T]}$-stopping time $R$ such that $S\leq
R\leq T$, we have $$
\dcb
\mathcal{G}^{(S,T]}_R&=&\{R=T\}\cap\mathcal{G}^*_{T}+\{R<T\}\cap\mathcal{G}_{R}.
\dce
$$
In particular, $\mathcal{G}^{(S,T]}_S=\mathcal{G}^{(S,T]}_0$ and
$\mathcal{G}^{(S,T]}_T=\mathcal{G}^*_{T}$.
\item[(3)]
For any $\mathbb{G}$-adapted process
$X$ such that $X_T\in\mathcal{G}^*_T$, $X^{(S,T]}$ is a
$\mathbb{G}^{(S,T]}$-adapted process. Conversely, for any
$\mathbb{G}^{(S,T]}$-adapted process $X'$, $X'^{(S,T]}$ is a
$\mathbb{G}$-adapted process. 
\item[(4)]  
For any $\mathbb{G}$-predictable process $K$, $K^{(S,T]}$ and
$\ind_{(S,T]}K$ define $\mathbb{G}^{(S,T]}$-predictable processes.
Conversely, for any $\mathbb{G}^{(S,T]}$-predictable process $K'$,
$K'^{(S,T]}$ and
$\ind_{(S,T]}K'$ define $\mathbb{G}$-predictable processes.  
\item[(5)]
For any $(\mathbb{Q},\mathbb{G})$ local martingale $X$ such
that $X_T\in\mathcal{G}^*_T$, $X^{(S,T]}$ is a
$(\mathbb{Q},\mathbb{G}^{(S,T]})$ local martingale. Conversely,
for any $(\mathbb{Q},\mathbb{G}^{(S,T]})$ local martingale $X'$,
$X'^{(S,T]}$ is a $(\mathbb{Q},\mathbb{G})$ local martingale.
\item[(6)]
For any $\mathbb{G}$-predictable process $K$, for
any $(\mathbb{Q},\mathbb{G})$ local martingale $X$ such that
$X_T\in\mathcal{G}^*_T$, the fact that $\ind_{(S,T]}K$ is
$X$-integrable in $\mathbb{G}$
implies that $\ind_{(S,T]}K$ is $X^{(S,T]}$-integrable in
$\mathbb{G}^{(S,T]}$. Conversely, for
any $\mathbb{G}^{(S,T]}$-predictable process $K'$, for any
$(\mathbb{Q},\mathbb{G}^{(S,T]})$ local martingale $X'$, the fact
that $\ind_{(S,T]}K'$ is $X'$-integrable in $\mathbb{G}^{(S,T]}$ implies that $\ind_{(S,T]}K'$ is
$X'^{(S,T]}$-integrable in $\mathbb{G}$.   
\dbe
\end{pro}

\subsection{$s\!\mathcal{H}$-measures}
\label{SHmeasure}

The following definition represents the notion of the "local immersion property" on a random interval under a change of probability measure.  

\begin{definition}
Let $S,T$ be $\mathbb{G}$-stopping times. A probability measure $\mathbb{Q}'$ defined on $\mathcal{G}_\infty$ is called an $s\!\mathcal{H}$-measure over the random time interval $(S,T]$ (with respect to $(\mathbb{Q},\mathbb{G})$), if $\mathbb{Q}'$ is equivalent to $\mathbb{Q}$ on $\mathcal{G}_{\infty}$, 
and if, for any $(\mathbb{Q},\mathbb{F})$ local martingale $X$, $X^{(S,T]}$ is a $(\mathbb{Q}',\mathbb{G}^{(S,T]})$ local martingale. 
\end{definition}

\brem
Note that, since obviously $X_T\in\mathcal{G}^*_T$, $X^{(S,T]}$ is $\mathbb{G}^{(S,T]}$-adapted (Proposition \ref{appendix} (3)). The notion of $s\!\mathcal{H}$-measure resembles to the immersion condition (also called $(\mathcal{H})$-hypothesis), but not exactly, because $\mathbb{Q}\neq \mathbb{Q}'$ and $\mathbb{F}\nsubseteq \mathbb{G}^{(S,T]}$. The notation $s\!\mathcal{H}$ refers to a "skewed" immersion condition. The $s\!\mathcal{H}$ measure condition was used to establish the $(\mathcal{H}')$ hypothesis in \cite{song}.
\erem

The notion of $s\!\mathcal{H}$-measure is closely linked with \mrt\ property.

\bethe\label{gSTMRT}
Suppose Assumption \ref{trois-hypos}. Then, for any $\mathbb{F}$-stopping time  $T$, for any $\mathbb{G}$-stopping time $S$ such that $S\geq \tau$, for any $s\!\mathcal{H}$-measure $\mathbb{Q}'$ over $(S,T]$, \mrt($\mathbb{Q}',\mathbb{G}^{(S,T]},W^{(S,T]}$) holds.
\ethe

\proof 
We suppose that the set $\{T>S\}$ is not empty, because otherwise, nothing is to be proved. The proof is presented in several steps.

$\blacktriangleright$ Let $\zeta$ be an $\mathcal{F}_{T}$-measurable bounded random variable, and $X$ be the martingale $X_t=\mathbb{Q}[\zeta|\mathcal{F}_t], 0\leq t < \infty$. The random time $T$ being an $\mathbb{F}$-stopping time, we have the identity $\zeta=X_{T\vee t}, t\geq 0$. By the \mrt($\mathbb{Q},\mathbb{F},W$) property, 
$X$ has the following representation in the filtration $\mathbb{F}$ under $\mathbb{Q}$:
\begin{equation}\label{integralrepresetation}
X_t=X_0+\int_0^tJ_s dW_s=X_0+\int_0^{t\wedge T}J_s dW_s,\ 0\leq t<\infty,
\end{equation}
where $J$ is a process in $\mathcal{I}(\mathbb{Q},\mathbb{F},W)$. 

$\blacktriangleright$ Since $W$ and $J{\stocint}W$ are $(\mathbb{Q},\mathbb{G})$ semimartingales, according to Lemma \ref{jfg}, $J$ is integrable with respect to $W$ in the sense of $(\mathbb{Q},\mathbb{G})$ semimartingale. This, together with Lemma \ref{equalintegral}, entails that the formula (\ref{integralrepresetation}) is also valid in the filtration $\mathbb{G}$. 

$\blacktriangleright$ We are in particular interested in a variant of the formula (\ref{integralrepresetation}):
\begin{equation}\label{STzeta}
X^{(S,T]}_t=\int_0^t\ind_{\{S<s\leq T\}}J_s dW^{(S,T]}_s,\ 0\leq t<\infty.
\end{equation}
We check straightforwardly that $J\ind_{(S,T]}$ is integrable with respect to $W^{(S,T]}$ in the sense of $(\mathbb{Q},\mathbb{G}^{(S,T]})$ semimartingale and the formula (\ref{STzeta}) is also valid in the filtration $\mathbb{G}^{(S,T]}$ under $\mathbb{Q}$. (cf. Proposition \ref{appendix} (6) and Lemma \ref{equalintegral}.)

$\blacktriangleright$ Let $\mathbb{Q}'$ be an $s\!\mathcal{H}$-measure on $(S,T]$. The two probabilities $\mathbb{Q}'$ and $\mathbb{Q}$ are equivalent on $\mathcal{G}_{\infty}$ and the processes $W^{(S,T]}$ and $X^{(S,T]}$ are $(\mathbb{Q}',\mathbb{G}^{(S,T]})$ local martingales. Recall that the (optional) brackets are the same under $\mathbb{Q}$ or under $\mathbb{Q}'$. We have$$
\ind_{(S,T]}J\transp [W,W\transp]J=[X^{(S,T]}].
$$
Note that $[X^{(S,T]}]$ is $\mathbb{G}^{(S,T]}$-adapted (cf. Proposition \ref{appendix} (3)). As $X$ is bounded, the bracket $[X^{(S,T]}]$ is $\mathbb{G}^{(S,T]}$-locally bounded, which entails $\ind_{(S,T]}J\in\mathcal{I}(\mathbb{Q}',\mathbb{G}^{(S,T]},W^{(S,T]})$.

From Lemma \ref{equalintegral}, we conclude that the formula (\ref{STzeta}) is also valid in the filtration $\mathbb{G}^{(S,T]}$ under $\mathbb{Q}'$.

$\blacktriangleright$ Consider now a bounded Borel function $g$ on $[0,\infty]$. Since $\tau\leq S$, we have $g(\tau)\in\mathcal{G}^*_S\subset\mathcal{G}^{(S,T]}_0$. Therefore, the process $g(\tau)\ind_{(S,T]}J$ is a $\mathbb{G}^{(S,T]}$-predictable process whilst $g(\tau)X^{(S,T]}$ is a $(\mathbb{Q}',\mathbb{G}^{(S,T]})$ martingale. Since $g$ is bounded, $g(\tau)\ind_{(S,T]}J\in\mathcal{I}(\mathbb{Q}',\mathbb{G}^{(S,T]},W^{(S,T]})$.  We can apply the formula (\ref{STzeta}) and write$:$ $$
g(\tau)X^{(S,T]}_t= \int_0^tg(\tau)\ind_{\{S<s\leq T\}}J_s dW^{(S,T]}_s,\ t\geq 0,
$$
valid in $\mathbb{G}^{(S,T]}$ under $\mathbb{Q}'$. 

$\blacktriangleright$ Note that $g(\tau)J\ind_{(S,T]}{\stocint} W^{(S,T]}$ is a bounded $(\mathbb{Q}',\mathbb{G}^{(S,T]})$ martingale. Also note that $X_t=X_T=\zeta$ for $t\geq T$ and $g(\tau)X_S\in\mathcal{G}^{(S,T]}_0$. We can write  $$
\dcb
&&g(\tau)X_{S}+\int_0^tg(\tau)\ind_{\{S<s\leq T\}}J_s dW^{(S,T]}_s\\
&=&\mathbb{Q}'[g(\tau)X_{S}+\int_0^\infty g(\tau)\ind_{\{S<s\leq T\}}J_s dW^{(S,T]}_s\ |\mathcal{G}^{(S,T]}_t]\\

&=&\mathbb{Q}'[g(\tau)\zeta\ |\mathcal{G}^{(S,T]}_t], t\geq 0.
\dce
$$ 
We deduce from this identity that $\mathbb{Q}'[g(\tau)\zeta\ |\mathcal{G}^{(S,T]}_0]=g(\tau)X_{S}$ and hence$$
\dcb
(g(\tau)\zeta)^\dag_{t}
=\int_0^tg(\tau)\ind_{\{S<s\leq T\}}J_s dW^{(S,T]}_s,\ t\geq 0,\\
\dce
$$
valid in $\mathbb{G}^{(S,T]}$ under $\mathbb{Q}'$, where we use the operator $\ii$ (as it is defined in Lemma \ref{boundedenough}) with respect to $(\mathbb{Q}',\mathbb{G}^{(S,T]})$. This identity shows that $(g(\tau)\zeta)^\dag$ belongs to $\mathcal{M}_{loc,0}(\mathbb{Q}',\mathbb{G}^{(S,T]},W^{(S,T]})$.

$\blacktriangleright$ Let $\mathcal{C}$ denote the class of all sets of the form $$
\{S<T\}\cap\{s<\tau\leq t\}\cap A+\{T\leq S\}\cap B,
$$ 
where $0\leq s<t$ and $A\in\mathcal{F}_T$ and $B\in\mathcal{G}^*_S$. Note that $\{S<T\},\{T\leq S\}\cap B\in\mathcal{G}^{(S,T]}_S=\mathcal{G}^{(S,T]}_0$ according to Proposition \ref{appendix} (2). This observation together with what we have proved previously shows that, for any $F\in\mathcal{C}$, the $(\mathbb{Q}',\mathbb{G}^{(S,T]})$ martingale $(\ind_F)^\dag$ belongs to $\mathcal{M}_{loc,0}(\mathbb{Q}',\mathbb{G}^{(S,T]},W^{(S,T]})$. On the other hand, the class $\mathcal{C}$ is a $\pi$-system and $\sigma(\mathcal{C})=\mathcal{G}^{(S,T]}_\infty$ because of Lemma \ref{g-star} and $\tau\leq S$. We apply Lemma \ref{boundedenough}, and we conclude with the property \mrt($\mathbb{Q}',\mathbb{G}^{(S,T]},W^{(S,T]}$). \ok

\

\subsection{From local solution to global solution}
\label{localglobal}

The next theorem shows how the local solutions are aggregated into a global one.

\bethe\label{mrt_after_default}
Suppose Assumption \ref{trois-hypos}. Assume moreover the following 

\textbf{$s\!\mathcal{H}$-measure condition covering $(\tau,\infty)$} : There exists a countable family of $\mathbb{G}$-stopping times $\{S_j,T_j: j\in \mathbb{N}\},$ such that
\ebe
\item
for any $j\in\mathbb{N}$, there exists an $s\!\mathcal{H}$-measure $\mathbb{Q}_j$ over the time interval $(S_j,T_j]$,

\item
$T_j$ are $\mathbb{F}$-stopping times,
\item
$S_j\geq \tau$ and
$(\tau,\infty)=\cup_{j\in\mathbb{N}}(S_j,T_j)$. 
\dbe
Then, \mrt($\mathbb{Q},\mathbb{G}^{(\tau,\infty]},\widetilde{W}^{(\tau,\infty]})$ holds. 
\ethe

Note that $\mathcal{G}^{(\tau,\infty]}_t=\mathcal{G}_{\tau\vee t}, t\geq 0$.

\proof The proof consists to establish firstly \mrt($\mathbb{Q},\mathbb{G}^{(S_j,T_j]},\widetilde{W}^{(S_j,T_j]})$ (the local solutions), and then to pass from local solutions to the global one. It is divided into several steps.

\textbf{1) Local solutions}

$\blacktriangleright$ Let $j\in\mathbb{N}$ be fixed and $\eta_t=\left.\frac{d\mathbb{Q}}{d\mathbb{Q}_j}\right|_{\mathcal{G}^{(S_j,T_j]}_t}, t\geq 0$. In this step, we prove that, for any $1\leq i\leq d$, $\cro{\eta,W^{(S_j,T_j]}_i}$ exists in $\mathbb{G}^{(S_j,T_j]}$ under $\mathbb{Q}_j$. To do this, it is enough to show that $(\eta W^{(S_j,T_j]}_{i})^\ast$ (where, for a process $\Xi$, $\Xi^*=\sup_{s\leq t}|\Xi_s|, t\geq 0$) is $\mathbb{G}^{(S_j,T_j]}$ locally $\mathbb{Q}_j$-integrable. Since $W^{(S_j,T_j]}_i$ is a $(\mathbb{Q},\mathbb{G}^{(S_j,T_j]})$ special semimartingale, there exists an increasing sequence of $\mathbb{G}^{(S_j,T_j]}$-stopping times $(R_n)$ converging to infinity $\mathbb{Q}$-almost surely (hence $\mathbb{Q}_j$ almost surely ) such that $(W^{(S_j,T_j]}_{i})^*_{R_n}$ is $\mathbb{Q}$-integrable for any $n\geq 1$. For $m\in\mathbb{N}$, set $$
U_m=\inf\{s: |(\eta W^{(S_j,T_j]}_{i})_s|>m \}.
$$ 
Then, $U_m$ tends to infinity $\mathbb{Q}_j$ almost surely and $$
|(\eta W^{(S_j,T_j]}_{i})^*_{U_m\wedge R_n}|\leq m+\eta_{U_m\wedge R_n}W^*_{U_m\wedge R_n}.
$$
The random variable on the right hand side is $\mathbb{Q}_j$-integrable, because $W^*_{U_m\wedge R_n}$ is $\mathbb{Q}$ integrable. This is what had to be proved.

$\blacktriangleright$
According to Theorem \ref{gSTMRT}, \mrt($\mathbb{Q}_j,\mathbb{G}^{(S_j,T_j]},W^{(S_j,T_j]})$ holds. By Lemma \ref{Girsanovrepresentation}, for $1\leq i\leq d$, \mrt($\mathbb{Q},\mathbb{G}^{(S_j,T_j]},W^{[\eta],(S_j,T_j]})$ holds, where $W^{[\eta],(S_j,T_j]}_i=W^{(S_j,T_j]}_i - \frac{1}{\eta_-}\stocint\cro{\eta,W^{(S_j,T_j]}_i}$ is a $(\mathbb{Q},\mathbb{G}^{(S_j,T_j]})$ local martingale. On the other hand, the process $$
\widetilde{W_i}^{(S_j,T_j]}=(W_i-\pv(W_i))^{(S_j,T_j]}=W^{(S_j,T_j]}_i-\pv(W_i)^{(S_j,T_j]}
$$ 
is also a $(\mathbb{Q},\mathbb{G}^{(S_j,T_j]})$ local martingale (cf. Proposition \ref{appendix} (5)). This implies that  $W^{[\eta],(S_j,T_j]}_i=\widetilde{W_i}^{(S_j,T_j]}$ and \mrt($\mathbb{Q},\mathbb{G}^{(S_j,T_j]},\widetilde{W}^{(S_j,T_j]})$ holds. 

\textbf{2) Global solution}

$\blacktriangleright$
Let us now prove the global solution \mrt($\mathbb{Q},\mathbb{G}^{(\tau,\infty]},\widetilde{W}^{(\tau,\infty]})$. Let $N$ be any element in $\mathcal{M}^\infty_{0}(\mathbb{Q},\mathbb{G}^{(\tau,\infty]})$ such that $N\widetilde{W}^{(\tau,\infty]}_i\in\mathcal{M}_{loc,0}(\mathbb{Q},\mathbb{G}^{(\tau,\infty]})$, $1\leq i\leq d$. We will prove that $N\equiv 0$, which, from Lemma \ref{caracMRT}, will achieve the proof of the theorem.

A "natural idea" would be to say that $N^{(S_j,T_j]}$ is a $(\mathbb{Q},\mathbb{G}^{(S_j,T_j]})$ local martingale, so that $N^{(S_j,T_j]}$ takes the form : $N^{(S_j,T_j]}=J\stocint \widetilde{W}^{(S_j,T_j]}$ because of \mrt($\mathbb{Q},\mathbb{G}^{(S_j,T_j]},\widetilde{W}^{(S_j,T_j]})$. This would entail $$
[N^{(S_j,T_j]}]=J\ind_{(S_j,T_j]}{\stocint}[N,\widetilde{W}]
\in
\mathcal{M}^\infty_{0}(\mathbb{Q},\mathbb{G}^{(\tau,\infty]}).
$$
Since $[N^{(S_j,T_j]}]$ is an increasing process, it would be null. We could then conclude $N\equiv 0$ by the covering condition.

However, this "natural idea" can not work because there is no guarantee that $N^{(S_j,T_j]}$ is a $(\mathbb{Q},\mathbb{G}^{(S_j,T_j]})$ local martingale. Proposition \ref{appendix} (5) is not applicable, if $\Delta_{S_j\vee T_j}N\notin\mathcal{G}^*_{S_j\vee T_j}$. That is why the proof of the global solution is concentrated on the study of $\Delta N$.

$\blacktriangleright$ 
We note that $N$ is a local martingale in the filtration $\mathbb{G}$. 
Indeed, $N_\tau=N_0=0$ because $\mathcal{G}^{(\tau,\infty]}_0=\mathcal{G}^{(\tau,\infty]}_\tau$ according to Proposition \ref{appendix} (2). This entails $N=N^{(\tau,\infty]}$. According to Proposition \ref{appendix} (5), $N$ is a $(\mathbb{Q},\mathbb{G})$ local martingale (as $N$ is bounded, it is in fact a true martingale). For the same reason $N\widetilde{W}^{(\tau,\infty]}_i$ is a $(\mathbb{Q},\mathbb{G})$ local martingale. Applying the integration by parts formula, we see that $
[N,\widetilde{W}_i]\in\mathcal{M}_{loc,0}(\mathbb{Q},\mathbb{G})
$.
Taking the stochastic integrals, we obtain that, for any process $J\in \mathcal{I}(\mathbb{Q},\mathbb{G},\widetilde{W})$,\begin{equation}\label{brac}
[N,J\stocint\widetilde{W}]\in\mathcal{M}_{loc,0}(\mathbb{Q},\mathbb{G}).
\end{equation}

$\blacktriangleright$ 
Let us study the jump process $\Delta N$ at predictable times. Let $T$ be a $\mathbb{G}$-predictable stopping time. We have $$
[{\ind_{[T]}}\stocint N,J\stocint\widetilde{W}]
=[ N,\ind_{[T]}J\stocint\widetilde{W}]
\in\mathcal{M}_{loc,0}(\mathbb{Q},\mathbb{G}).
$$
Consider the $(\mathbb{Q},\mathbb{G})$ local martingale $X=\ind_{[T]}{\stocint} N$. We have
$
X=\ind_{\{\tau<T\}}\Delta_{T}N\ind_{[T,\infty)}
$
and, for any $j\in\mathbb{N}$, 
$
\Delta_{T_j}X^{T_j}=\ind_{\{T=T_j\}}\Delta_{T}N
$.

Let $j\in\mathbb{N}$ be fixed. Set $\varkappa=\ind_{\{T=T_j\}}\Delta_{T}N$ and $\varkappa_-=\mathbb{Q}[\varkappa|\mathcal{G}_{T-}]$. Using the fact that $T$ is $\mathbb{G}$-predictable, by a direct computation, we obtain 
$$
(\varkappa\ind_{[T_j,\infty)})^{(p)}=\ind_{\{T\leq T_j\}}\varkappa_-\ind_{[T,\infty)},
$$
where $(\varkappa\ind_{[T_j,\infty)})^{(p)}$ denotes the $(\mathbb{Q},\mathbb{G})$ predictable dual projection of the jump process $\varkappa\ind_{[T_j,\infty)}$. 

Set $X'=\varkappa\ind_{[T_j,\infty)}-(\varkappa\ind_{[T_j,\infty)})^{(p)}$ and $$
X''=X^{T_j}-X'=
(\ind_{\{T< T_j\}}\Delta_{T}N +\ind_{\{T\leq T_j\}}\varkappa_-)\ind_{[T,\infty)}.
$$ 
Since $X''=X''^{T_j}$ and $T$ is $\mathbb{G}$-predictable,$$
X''_{S_j\vee T_j}=X''_{T_j}=\ind_{\{T< T_j\}}\Delta_{T}N +\ind_{\{T\leq T_j\}}\varkappa_-\in\mathcal{G}_{T_j-}\subset \mathcal{G}_{S_j\vee T_j-}
\subset \mathcal{G}^*_{S_j\vee T_j}.
$$
Proposition \ref{appendix} (5) implies that $X''^{(S_j,T_j]}\in\mathcal{M}_{loc,0}(\mathbb{Q},\mathbb{G}^{(S_j,T_j]})$. By \mrt($\mathbb{Q},\mathbb{G}^{(S_j,T_j]},\widetilde{W}^{(S_j,T_j]})$, there exists a $d$-dimensional process $J\in\mathcal{I}(\mathbb{Q},\mathbb{G}^{(S_j,T_j]},\widetilde{W}^{(S_j,T_j]})$ such that $X''^{(S_j,T_j]} = J\stocint \widetilde{W}^{(S_j,T_j]}$ in $(\mathbb{Q},\mathbb{G}^{(S_j,T_j]})$. Applying Proposition \ref{appendix} (4) and (6) and Lemma \ref{equalintegral}, we have $J\ind_{(S_j,T_j]}\in\mathcal{I}(\mathbb{Q},\mathbb{G},\widetilde{W})$ and $$
X''^{(S_j,T_j]} = J\ind_{(S_j,T_j]}{\stocint} \widetilde{W}^{(S_j,T_j]}=J\ind_{(S_j,T_j]}{\stocint} \widetilde{W}
$$ 
in $(\mathbb{Q},\mathbb{G})$. From this relation, in computing the jump at $T$, we deduce$$
(\ind_{\{T< T_j\}}\Delta_{T}N +\ind_{\{T\leq T_j\}}\varkappa_-) \ind_{\{S_j<T\leq T_j\}}=J_T \ind_{\{S_j<T\leq T_j\}}\Delta_T\widetilde{W}.
$$
We now compute the bracket between $N$ and ${J\ind_{(S_j,T_j]}}\stocint\widetilde{W}$ (see (\ref{brac})):$$
\dcb
[{\ind_{[T]}}\stocint N,J\ind_{(S_j,T_j]}{\stocint}\widetilde{W}]

&=&(\ind_{\{T< T_j\}}(\Delta_{T}N)^2 +\Delta_TN\ind_{\{T\leq T_j\}}\varkappa_-)\ind_{\{S_j<T\leq T_j\}}\ind_{[T,\infty)}.\\
\dce
$$
We know that $[{\ind_{[T]}}\stocint N,J\ind_{(S_j,T_j]}{\stocint}\widetilde{W}]$ is a $(\mathbb{Q},\mathbb{G})$ local martingale. We check that $$
\Delta_TN\ind_{\{T\leq T_j\}}\varkappa_-\ind_{\{S_j<T\leq T_j\}}\ind_{[T,\infty)}
$$ 
also is a $(\mathbb{Q},\mathbb{G})$ local martingale. We conclude therefore $\Delta_{T}N\ind_{\{S_j<T<T_j\}}=0$.

This nullity being true for any $j\in\mathbb{N}$, we can apply the covering condition and we conclude $\Delta_{T}N=\Delta_{T}N\ind_{\{\tau<T<\infty\}}=0$.

$\blacktriangleright$
We have proved that $N$ has no jumps at $\mathbb{G}$-predictable times. For $j\in\mathbb{N}$, we introduce the process $$
N'=\Delta_{S_j\vee T_j}N\ind_{[S_j\vee T_j,\infty)} -  (\Delta_{ S_j\vee T_j}N\ind_{[S_j\vee T_j,\infty)})^{(p)}.
$$ 
The process $(\Delta_{S_j\vee T_j}N\ind_{[S_j\vee T_j,\infty)})^{(p)}$ is continuous, because $N$ has no jump at predictable times. Set $N''=N-N'$ and compute the jump of $N''$ at $S_j\vee T_j$: $$
\Delta_{S_j\vee T_j}N''=\Delta_{S_j\vee T_j}N-\Delta_{S_j\vee T_j}N'=0.
$$
This nullity entails that $(N'')_{S_j\vee T_j}\in\mathcal{G}^*_{S_j\vee T_j}$, so that $N''^{(S_j,T_j]}\in\mathcal{M}_{loc,0}(\mathbb{Q},\mathbb{G}^{(S_j,T_j]})$ (cf. Proposition \ref{appendix} (5)). Because of the property \mrt($\mathbb{Q},\mathbb{G}^{(S_j,T_j]},\widetilde{W}^{(S_j,T_j]})$, there exists a $d$-dimensional process $J\in\mathcal{I}(\mathbb{Q},\mathbb{G}^{(S_j,T_j]},\widetilde{W}^{(S_j,T_j]})$ such that $N''^{(S_j,T_j]} = J\stocint \widetilde{W}^{(S_j,T_j]}$. Applying Proposition \ref{appendix} (6), Lemma \ref{equalintegral}, we have also $J\ind_{(S_j,T_j]}\in\mathcal{I}(\mathbb{Q},\mathbb{G},\widetilde{W})$ and $N''^{(S_j,T_j]} = J\ind_{(S_j,T_j]}{\stocint} \widetilde{W}$ in the sense of $(\mathbb{Q},\mathbb{G})$. On the other hand, we check immediately $[N''^{(S_j,T_j]},N'^{(S_j,T_j]}]=0$. These facts enable us to write
$$
[N''^{(S_j,T_j]},N''^{(S_j,T_j]}]=[N,J\ind_{(S_j,T_j]}{\stocint}\widetilde{W}_i]\in\mathcal{M}_{loc,0}(\mathbb{Q},\mathbb{G}).
$$
This relation is possible only if $N''^{(S_j,T_j]}\equiv 0$, which yields $N^{(S_j,T_j]}=N'^{(S_j,T_j]}$. It follows that $\ind_{\{S_j<t<T_j\}}\Delta_tN=0$ and $N$ has bounded variation on $(S_j,T_j]$.

$\blacktriangleright$ 
Now, by covering condition $(\tau,\infty)=\cup_{j\in\mathbb{N}}(S_j,T_j)$, we conclude that $N$ is a continuous local martingale with finite variation. It is therefore a constant, i.e., it is null. \ok

\

\section{\mrt\ on $\mathbb{R}_+$ and an equality between $\mathcal{G}_{\tau-}$ and $\mathcal{G}_{\tau}$}\label{sectionG=G}

The \mrt\ property in the progressively enlarged filtration is closely linked with the $\sigma$-algebra equality $\{0<\tau<\infty\}\cap \mathcal{G}_{\tau}=\{0<\tau<\infty\}\cap \mathcal{G}_{\tau-}$. The study of the gap between $\mathcal{G}_{\tau-}$ and $\mathcal{G}_{\tau}$ has long been considered because of its importance in progressive enlargement of filtration (see for example \cite{BEKSY,J,song-splitting}). Our discussions below give complementary information to this problem. 

Putting together Theorem \ref{companion} and Theorem \ref{mrt_after_default}, we obtain immediately

\bethe\label{put-together}
Suppose Assumption \ref{trois-hypos} and $s\!\mathcal{H}$-measure condition covering $(\tau,\infty)$. Then, \mrt$(\mathbb{Q},\mathbb{G},(\widetilde{W},L))$  and $\ind_{\{0<\tau<\infty\}}W_\tau\in\mathcal{G}_{\tau-}$ hold if and only if $\{0<\tau<\infty\}\cap\mathcal{G}_\tau=\{0<\tau<\infty\}\cap\mathcal{G}_{\tau-}$.
\ethe

We see the particular role played by the equality $\{0<\tau<\infty\}\cap\mathcal{G}_{\tau-}=\{0<\tau<\infty\}\cap\mathcal{G}_{\tau}$. In this section we show how this equality can be studied by $s\!\mathcal{H}$-measure condition.

\bethe\label{G=G}
Suppose Assumption \ref{trois-hypos} and the following 

\textbf{$s\!\mathcal{H}$-measure condition covering $(0,\infty)$}:
 There exists a countable family of $\mathbb{G}$-stopping times $\{S_j,T_j: j\in \mathbb{N}\},$ such that
\ebe
\item
for any $j\in\mathbb{N}$, there exists an $s\!\mathcal{H}$-measure $\mathbb{Q}_j$ over the time interval $(S_j,T_j]$.

\item
$T_j$ are $\mathbb{F}$-stopping times.
\item
$(0,\infty)=\cup_{i\in\mathbb{N}}(S_j,T_j)$. 
\dbe
Suppose that $\tau$ avoids the $\mathbb{F}$ stopping times on $(0,\infty)$. Then, $\{0<\tau<\infty\}\cap\mathcal{G}_{\tau-}=\{0<\tau<\infty\}\cap\mathcal{G}_{\tau}$.
\ethe

\proof Fix $j\in\mathbb{N}$. Let $\zeta$ be an $\mathcal{F}_{T_j}$-measurable bounded random variable, and $X$ be the martingale $X_t=\mathbb{Q}[\zeta|\mathcal{F}_t], 0\leq t < \infty$. We note the identity $\zeta=X_{T_j\vee t}$ for all $t\geq 0$. By definition of $s\!\mathcal{H}$ measure, $X^{(S_j,T_j]}$ is  a $(\mathbb{Q}_j,\mathbb{G}^{(S_j,T_j]})$ uniformly integrable martingale. Let $R_j=S_j\vee(\tau\wedge T_j)$. Note that $\{S_j<R_j<T_j\}$ is equivalent to $\{S_j<\tau<T_j\}$ (and in particular $R_j=\tau$). Let $g$ be a bounded Borel function. Applying Lemma \ref{gSTstoppingtimes} and Proposition \ref{appendix} (2), we can write
$$
\dcb
&&\ind_{\{S_j<\tau<T_j\}}\mathbb{Q}_j[g(\tau)\zeta|\mathcal{G}_{\tau}]

=
g(\tau)\ind_{\{S_j<R_j<T_j\}}\mathbb{Q}_j[X^{S_j\vee T_j}_{S_j\vee T_j}|\mathcal{G}^{(S_j,T_j]}_{R_j}]\\

&=&g(\tau)\ind_{\{S_j<\tau<T_j\}}X_{\tau}

=
g(\tau)\ind_{\{S_j<\tau<T_j\}}X_{\tau-}\in \mathcal{G}_{\tau-},\\ 
&&\mbox{(because $\tau$ avoids the $\mathbb{F}$ stopping times and $\mathbb{Q}_j$ is equivalent to $\mathbb{Q}$).}
\dce
$$
Since $\{S_j<\tau<T_j\}\cap\mathcal{G}_\tau\subset\{S_j<\tau<T_j\}\cap(\sigma(\tau)\vee\mathcal{F}_{T_j})$, the above relation yields $\{S_j<\tau<T_j\}\cap\mathcal{G}_{\tau}=\{S_j<\tau<T_j\}\cap\mathcal{G}_{\tau-}$ under $\mathbb{Q}_j$ (and hence, under $\mathbb{Q}$). 

Now, for any $A\in\mathcal{G}_\tau$, noting that $(0,\infty)=\cup_{i\in\mathbb{N}}(S_j,T_j)$, we can write under the probability $\mathbb{Q}$ $$
\dcb
\{0<\tau<\infty\}\cap A
&=&\cup_{j\in\mathbb{N}}\{S_j<\tau<T_j\}\cap A\\
&=&\cup_{j\in\mathbb{N}}\{S_j<\tau<T_j\}\cap B_j\ \mbox{ for some $B_j\in\mathcal{G}_{\tau-}$}\\
&=&\cup_{j\in\mathbb{N}}\{S_j<\tau\leq T_j, \tau<\infty\}\cap B_j\ 
\mbox{ because $\tau$ avoids the $\mathbb{F}$ stopping times}\\
&\in&\{0<\tau<\infty\}\cap\mathcal{G}_{\tau-}.
\dce
$$
We have proved $\{0<\tau<\infty\}\cap\mathcal{G}_{\tau}\subset \{0<\tau<\infty\}\cap\mathcal{G}_{\tau-}$ under $\mathbb{Q}$. The inverse inclusion being an evidence, the theorem is proved. \ok

We end this section by the following relation between $s\!\mathcal{H}$-measure condition covering $(0,\infty)$ and $s\!\mathcal{H}$-measure condition covering $(\tau,\infty)$. 

\bl\label{reducing}
If the family $\{S_j,T_j: j\in\mathbb{N}\}$ of $\mathbb{G}$ stopping times satisfies the $s\!\mathcal{H}$-measure condition covering $(0,\infty)$, the family $\{(S_j\vee\tau)\wedge (S_j\vee T_j),\ T_j: j\in\mathbb{N}\}$ satisfies the $s\!\mathcal{H}$-measure condition covering $(\tau,\infty)$.
\el

\proof Consider four $\mathbb{G}$ stopping times $S,T,U,V$. Suppose $S\leq U\leq S\vee T, S\leq V\leq S\vee T$. By Proposition \ref{appendix}, Lemma \ref{gSTstoppingtimes}, for any $\mathbb{G}$ stopping time $R$ such that $U\leq R\leq U\vee V$, one has $\mathcal{G}^{(U,V]}_{R}\subset\mathcal{G}^{(S,T]}_{R}$.

Let $\mathbb{Q}'$ be an $s\!\mathcal{H}$-measure on $(S,T]$. For any $(\mathbb{Q},\mathbb{F})$ uniformly integrable martingale $X$, $X^{(S,T]}$ is a $(\mathbb{Q}',\mathbb{G}^{(S,S\vee T]})$ uniformly integrable martingale, and consequently, $X^{(U,V]}=(X^{(S,T]})^{(U,V]}$ is a $(\mathbb{Q}',\mathbb{G}^{(S,S\vee T]})$ uniformly integrable martingale. For any $\mathbb{G}$ stopping times $R,R'$ such that $U\leq R\leq R'\leq U\vee V$, for $A\in\mathcal{G}^{(U,U\vee V]}_R$, we have $A\in\mathcal{G}^{(S,S\vee T]}_R$ and therefore $$
\dcb
\mathbb{Q}'[\ind_AX^{(U,V]}_{R'}]
=\mathbb{Q}'[\ind_A(X^{(S,T]})^{(U,V]}_{R'}]
=\mathbb{Q}'[\ind_A(X^{(S,T]})^{(U,V]}_{R}]
=\mathbb{Q}'[\ind_AX^{(U,V]}_{R}].
\dce
$$
Since $X^{(U,V]}$ is $\mathbb{G}^{(U,U\vee V]}$ adapted, $X^{(U,V]}$ is a $(\mathbb{Q}',\mathbb{G}^{(U,U\vee V]})$ uniformly integrable martingale. This is enough to conclude the lemma. \ok

\

\section{Examples}
\label{examples}

In this section we assume the setting of subsection \ref{framework}. We show how our method applies, in the situations already studied in the literature, to provide a uniform way to prove various classical results. In all this section, Assumption \ref{trois-hypos} (i) is in force.

\subsection{The case of the immersion condition}\label{sect_H}

Suppose the immersion condition (\cite{BJR,BY,Kusuoka}), i.e., any $(\mathbb{Q},\mathbb{F})$ local martingale is a $(\mathbb{Q},\mathbb{G})$ local martingale (in particular $W=\widetilde{W}$). In this case, if we take $T=\infty$ and $S=0$, the random variable $T$ is an $\mathbb{F}$ stopping time and the interval $(S,T]$ covers $(0,\infty)$, and the probability measure $\mathbb{Q}$ is clearly an $s\!\mathcal{H}$-measure on $(S,T]$. Hence, according to Lemma \ref{reducing} and Theorem \ref{put-together}, the properties \mrt$(\mathbb{Q},\mathbb{G},(W,L))$ and $\ind_{\{0<\tau<\infty\}}W_{\tau}\in\mathcal{F}_{\tau-}$ hold, whenever $\{0<\tau<\infty\}\cap\mathcal{G}_\tau=\{0<\tau<\infty\}\cap\mathcal{G}_{\tau-}$.

Let us show that the immersion condition together with $W_{\tau}\in\mathcal{F}_{\tau-}$ implies $\mathcal{G}_{\tau-}=\mathcal{G}_{\tau}$. Let $\zeta$ be a bounded $\mathcal{F}_\infty$ measurable random variable. Let $\zeta_t=\mathbb{Q}[\zeta|\mathcal{F}_t], t\geq 0$. Since $W_{\tau}\in\mathcal{F}_{\tau-}$, by the property \mrt$(\mathbb{Q},\mathbb{F},W)$, we check that $\Delta_{\tau}\zeta\in\mathcal{F}_{\tau-}$ so that $\zeta_{\tau}\in\mathcal{F}_{\tau-}=\mathcal{G}_{\tau-}$ (see the proof of Theorem \ref{companion}). Thanks to the immersion condition, if $g$ is a bounded Borel function on $[0,\infty]$, $$
\mathbb{Q}[g(\tau)\zeta|\mathcal{G}_\tau]=g(\tau)\mathbb{Q}[\zeta|\mathcal{G}_\tau]
=g(\tau)\zeta_\tau \in \mathcal{G}_{\tau-}.
$$
Applying the monotone class theorem, we obtain $\mathcal{G}_{\tau}=\mathcal{G}_{\tau-}$.

\bethe\label{H_hypothesis}
Suppose Assumption \ref{trois-hypos} (i). Then the following two conditions are equivalent
\ebe
\item[(i)]
Immersion condition and $W_{\tau}\in\mathcal{F}_{\tau-}$.
\item[(ii)]
\mrt($\mathbb{Q},\mathbb{G},(W,L))$ and $\mathcal{G}_\tau=\mathcal{G}_{\tau-}$.
\dbe
\ethe

\brem
\ebe
\item[a.]
Theorem \ref{H_hypothesis} implies the results given in \cite[Theorem 2.3, $N=1$]{Kusuoka}, because, for $W$ a Brownian motion, $W_\tau\in\mathcal{F}_{\tau-}$. Applying Theorem \ref{H_hypothesis} and Lemma \ref{Girsanovrepresentation}, we can prove also \cite[Theorem 3.2, $N=1$]{Kusuoka}.
\item[b.]
The condition $W_{\tau}\in\mathcal{F}_{\tau-}$ holds obviously, if $\tau$ avoids the $\mathbb{F}$ stopping times.
\dbe
\erem

\subsection{The case of a honest time}\label{honest}

In this section, we suppose that $\tau$ is an $\mathbb{F}$ honest time. Honest times have been fully studied in the past (cf. \cite{barlow,J,JY2}). Let us reconsider this case with the results obtained in this paper. 

To simplify the computations, we assume

\noindent\textbf{Hy}(C) All $(\mathbb{Q},\mathbb{F})$ local martingales are continuous.

We know that, when $\tau$ is an honest time, $(\mathcal{H}')$ hypothesis holds. Under \textbf{Hy}(C), $\hat{A}-A\equiv 0$ (see Section \ref{framework} for notations) and the drift operator $\pv(X)$ on $(\tau,\infty)$ is given by$$
\ind_{(\tau,\infty)}{\stocint} \pv(X) = \ind_{(\tau,\infty)}\frac{1}{1-Z}{\stocint} \cro{M,X}
$$
which is continuous. 

\bl
Suppose Assumption \ref{trois-hypos} (i) and \textbf{Hy}(C). Then, we have the $s\!\mathcal{H}$-measure condition covering $(\tau,\infty)$.
\el

\proof For any $n\in\mathbb{N}^*$, for any $a>0$, set $$
T_{a,n}=\inf\{t\geq a: \int_a^t\frac{1}{(1-Z_{s-})^2} d\cro{M}_s>n\}\wedge n.
$$
The random variables $T_{a,n}$ are $\mathbb{F}$ stopping times. Set $S_{a}=\tau\vee a$. Since $Z_{t-}<1, Z_t<1,$ for $t>\tau$, we have $\cup_{a\in\mathtt{Q},a>0,n\in\mathbb{N}^*}(S_a,T_{a,n})=(\tau,\infty)$. For fixed $a>0, n\in\mathbb{N}^*$, we introduce the process$$
\eta = \mathcal{E}\Big(\ind_{(S_a,T_{a,n}]}\frac{1}{1-Z_{-}}{\stocint}\widetilde{M}\Big)
$$
which is a positive continuous $(\mathbb{Q},\mathbb{G})$ martingale. We now show that the probability measure $$
\mathbb{Q}^{a,n}=\eta_{S_a\vee T_{a,n}}\cdot\mathbb{Q}
$$
is an $s\!\mathcal{H}$-measure on the random interval $(S_a,T_{a,n}]$.

Let $X$ be a $(\mathbb{Q},\mathbb{F})$ local martingale. Then, $\widetilde{X}$ is a $(\mathbb{Q},\mathbb{G})$ local martingale. By Girsanov's theorem, the process $\widetilde{X}^{[\eta]}$ (cf. Lemma \ref{Girsanovrepresentation}) is a $(\mathbb{Q}^{a,n},\mathbb{G})$ local martingale. Let us compute $\widetilde{X}^{[\eta]}$ on the random interval $(S_a,T_{a,n}]$. Thanks to \textbf{Hy}(C), $\cro{\widetilde{X},\widetilde{M}}=\cro{X,M}$ and $\cro{\eta,\widetilde{X}}=\cro{\eta,X}$. 
$$
\dcb
&&\widetilde{X}^{[\eta]}_{\tau\vee t}-\widetilde{X}^{[\eta]}_{\tau\vee a}\\
&=&
X_{\tau\vee t}-X_{\tau\vee a}+\int_{a}^{t}\ind_{\{\tau< s\}}\frac{1}{1-Z_{s-}}d\cro{M,X}_s
-\int_{a}^{t}\ind_{\{\tau< s\}}\frac{1}{1-Z_{s-}}\ind_{\{S_a<s\leq T_{a,n}\}}d\cro{\widetilde{M},X}_s.\\
\\

&=&X_{\tau\vee t}-X_{\tau\vee a},\ \mbox{ if $a\leq t\leq T_n$}.
\dce
$$ 
It follows that $X^{\tau\vee a\vee T_{a,n}}_{\tau\vee a\vee t}-X_{\tau\vee a}, t\geq 0,$ is a $(\mathbb{Q}^{a,n},\mathbb{G})$ local martingale. Since $
X_{\tau\vee a\vee T_{a,n}}\in\mathcal{F}_{\tau\vee a\vee T_{a,n}}\subset \mathcal{G}^*_{\tau\vee a\vee T_{a,n}},
$
we conclude that $X^{(S_a,T_{a,n}]}$ is a $(\mathbb{Q}^{a,n},\mathbb{G}^{(S_a,T_{a,n}]})$ local martingale (cf. Proposition \ref{appendix} (5)). This proves that the probability $\mathbb{Q}^{a,n}$ is an $s\!\mathcal{H}$-measure on $(S_a,T_{a,n}]$. \ok

Now we apply  Theorem \ref{beforeTau} and Theorem \ref{mrt_after_default}. Note that, according to \cite[Proposition(5.3)]{J}, for any $\mathbb{G}$-predictable process $J$, there exist $\mathbb{F}$-predictable processes $J',J''$ such that$$
J\ind_{(0,\infty)}=J'\ind_{(0,\tau]}+J''\ind_{(\tau,\infty)}.
$$
Recall the process $L$ in Remark \ref{Ldef}.

\bethe
Suppose Assumption \ref{trois-hypos} (i) and \textbf{Hy}(C). For any bounded $(\mathbb{Q},\mathbb{G})$-martingale $X$, there exist $\mathbb{F}$-predictable processes $J', J'', K$ and a bounded $\xi\in\mathcal{G}_\tau$ such that $\mathbb{Q}[\xi|\mathcal{G}_{\tau-}]=0$, $J'\ind_{(0,\tau]}+J''\ind_{(\tau,\infty)}\in\mathcal{I}(\mathbb{Q},\mathbb{G},\widetilde{W})$ and, for $t\geq 0$, 
\begin{equation}\label{for-honest}
X_t=X_0+\int_0^t\left(J'\ind_{\{s\leq\tau\}}+J''\ind_{\{\tau<s\}}\right)d\widetilde{W}_s
+\int_0^tK_s\ind_{\{0<s\leq\tau\}}dL_s+\ind_{\{\tau>0\}}\xi H_t.
\end{equation}
If, in addition,  $\{0<\tau<\infty\}\cap\mathcal{G}_{\tau-}=\{0<\tau<\infty\}\cap\mathcal{G}_{\tau}$, the property \mrt($\mathbb{Q},\mathbb{G},(\widetilde{W},L))$ holds.
\ethe

We end this section by a remark on Brownian filtrations (in the sense of \cite{BEKSY}).

\bethe\label{black-scholes-completion}
Suppose Assumption \ref{trois-hypos} (i). Suppose that $\mathbb{F}$ is a Brownian filtration. Then there exists a bounded $(\mathbb{Q},\mathbb{G})$ martingale $\nu$ such that \mrt($\mathbb{Q},\mathbb{G},(\widetilde{W},L,\nu))$ holds.
\ethe

\proof According to \cite{BEKSY}, since $\mathbb{F}$ is a Brownian filtration and $\tau$ is $\mathbb{F}$ honest, there exists a random event $A\in\mathcal{G}_\tau$ such that $\mathcal{G}_\tau=\mathcal{G}_{\tau-}\vee\sigma(A)$. This means that, for any $\xi\in \mathcal{G}_\tau$, there exist $\xi', \xi''\in \mathcal{G}_{\tau-}$ such that $\xi=\xi'\ind_A+\xi''\ind_{A^c}$. In particular, if $\mathbb{Q}[\xi|\mathcal{G}_{\tau-}]=0$, i.e., $0 = \xi' p+\xi''(1-p)$, where $p=\mathbb{Q}[A|\mathcal{G}_{\tau-}]$, we have
$$
\dcb
\xi
&=&(-\ind_{\{p>0\}}\xi''\frac{1}{p}+\ind_{\{p=0\}}\xi'\frac{1}{1-p})((1-p)\ind_A-p\ind_{A^c}).\\
\dce
$$
Let $\nu=((1-p)\ind_A-p\ind_{A^c})H$ and $F$ be an $\mathbb{F}$-predictable process such that $F_\tau=(-\ind_{\{p>0\}}\xi''\frac{1}{p}+\ind_{\{p=0\}}\xi'\frac{1}{1-p})$. Then, $\nu$ is a bounded $(\mathbb{Q},\mathbb{G})$ martingale and $F$ is $\nu$-integrable. For bounded $(\mathbb{Q},\mathbb{G})$-martingale $X$, the formula (\ref{for-honest}) now becomes
$$
X_t=X_0+\int_0^t\left(J'\ind_{\{s\leq\tau\}}+J''\ind_{\{\tau<s\}}\right)d\widetilde{W}_s\\
+\int_0^tK_s\ind_{\{0<s\leq\tau\}}dL_s+\int_0^tF_s d\nu_s.
$$
This proves the theorem. \ok

\subsection{The case of density hypothesis}

In this subsection we work under 

\bassumption\label{density_hypo}
\textbf{Density Hypothesis.} We assume that, for any $t\in\mathbb{R}_+$, there exists a strictly positive $\mathcal{B}[0,\infty]\otimes\mathcal{F}_t$ measurable function $\alpha_t(\theta,\omega),(\theta,\omega)\in[0,\infty]\times\Omega$, which gives the conditional law $$
\mathbb{Q}[\tau \in \mathtt{A} |\mathcal{F}_t] = \int_\mathtt{A}\alpha_t(\theta)\mu(d\theta),\ t\geq 0, \mathtt{A}\in\mathcal{B}[0,\infty],
$$ 
where $\alpha_t(\theta)$ denotes the application $\alpha_t(\theta,\cdot)$ and $\mu$ is a diffuse probability measure on $\mathbb{R}_+$. We assume that the trajectory $t\rightarrow \alpha_t(\theta,\omega)$ is càdlàg. 
\eassumption

For any $n\geq 1$, let $\mathbb{Q}_n' = \alpha_n(\tau)^{-1}\cdot\mathbb{Q}$.
We check that $
\mathbb{Q}_n'[\tau\in d\theta|\mathcal{F}_n]=
\mu(d\theta).
$
This means that, under $\mathbb{Q}_n'$, $\tau$ is independent of $\mathcal{F}_n$. For any $(\mathbb{Q}_n',\mathbb{F})$-local martingale $Y$, the process $Y^{n}$ ($Y$ stopped at $n$) will be a $(\mathbb{Q}_n',\mathbb{G})$-local martingale, and therefore $Y^{(0,n]}$ is a $(\mathbb{Q}_n',\mathbb{G}^{(0,n]})$-local martingale. Since$$
\mathbb{Q}[\frac{1}{\alpha_n(\tau)}|\mathcal{F}_n]= 1,
$$
we have $\mathbb{Q}_n'|_{\mathcal{F}_n}=\mathbb{Q}|_{\mathcal{F}_n}$. This yields that, for any $(\mathbb{Q},\mathbb{F})$-local martingale $X$, $X^n$ is a $(\mathbb{Q}_n',\mathbb{F})$-local martingale, and $X^{(0,n]}$ is a $(\mathbb{Q}_n',\mathbb{G}^{(0,n]})$-local martingale. We have proved that $\mathbb{Q}_n'$ is an $s\!\mathcal{H}$-measure on $(0,n]$. We note that the integers $n$ are $\mathbb{F}$ stopping times and the intervals $(0,n],n\geq 1,$ cover $(0,\infty)$. Note also that, since $\mu$ is diffuse, $\tau$ avoids the $\mathbb{F}$ stopping times.

Applying Theorem \ref{G=G}, Lemma \ref{reducing}, Theorem \ref{put-together}, we obtain:

\bethe\label{DH_result}
Suppose Assumption \ref{trois-hypos} (i) and density hypothesis \ref{density_hypo}. Then, $\{0<\tau<\infty\}\cap\mathcal{G}_{\tau-}=\{0<\tau<\infty\}\cap\mathcal{G}_{\tau}$ and \mrt($\mathbb{Q},\mathbb{G},(\widetilde{W},L))$ holds.
\ethe 

\brem
The property \mrt\ under density hypothesis \ref{density_hypo} has been studied in \cite[Theorem 2.1]{JC2} when $W$ is continuous, making use of Itô's computations.
\erem

\subsection{The case of Cox measure and the related ones}

In this subsection, we consider a probability space $(\Omega,\mathcal{A},\mathbb{P})$ equipped with a filtration $\mathbb{F}$ of sub-$\sigma$-algebras in $\mathcal{A}$. We consider the product measurable space $([0,\infty]\times \Omega,\mathcal{B}[0,\infty]\otimes\mathcal{F}_\infty)$. As usual, we consider $\mathbb{F}$ as a filtration on the product space and $\mathbb{P}$ as a probability measure defined on $\mathcal{F}_\infty$ considered as a sub-$\sigma$-algebra of $\mathcal{B}[0,\infty]\otimes\mathcal{F}_\infty$  (see Section \ref{JSmodel}). Consider the projection map : $\tau(s,\omega)=s$ for $(s,\omega)\in [0,\infty]\times\Omega$. Let $\Lambda$ be a continuous increasing $\mathbb{F}$-adapted process such that $\Lambda_0=0,\Lambda_\infty=\infty$. The Cox measure $\nu^\Lambda$ on the product space $[0,\infty]\times\Omega$ associated with $\Lambda$ is defined by the relation
\begin{equation}\label{cox}
\nu^\Lambda[A \cap \{s<\tau\leq
t\}]=\mathbb{P}[\mathbb{I}_{A}\int_s^t e^{-\Lambda_v}d\Lambda_v],\
A\in\mathcal{F}_\infty,\ 0<s<t<\infty.
\end{equation} 
Consider the progressively enlarged filtration $\mathbb{G}$ on the product space $[0,\infty]\times\Omega$ from $\mathbb{F}$ with $\tau$. It is well know that, under the Cox measure, the immersion condition holds (cf. \cite{BJR}). It is also easy to check that $\nu^\Lambda[\tau=T]=0$ for any $\mathbb{F}$-stopping time $T$, consequence of the continuity of the process $\Lambda$. This last property implies $W_{\tau}\in\mathcal{F}_{\tau-}$. Theorem \ref{H_hypothesis} is applicable. We have the property \mrt($\nu^\Lambda,\mathbb{G},(W,H-\frac{1}{1-e^{-\Lambda}}\stocint \Lambda))$ and $\mathcal{G}_{\tau-}=\mathcal{G}_{\tau}$.

Now, if a probability measure $\mathbb{Q}$ on the product space is absolutely continuous with respect to the Cox measure, we apply Lemma \ref{Girsanovrepresentation} to obtain 

\bethe\label{Cox_result}
Suppose Assumption \ref{trois-hypos} (i). If the probability measure $\mathbb{Q}$ is absolutely continuous with respect to the Cox measure $\nu^\Lambda$, we have the properties $W_{\tau}\in\mathcal{F}_{\tau-}$, $\mathcal{G}_{\tau-}=\mathcal{G}_{\tau}$ and  \mrt($\mathbb{Q},\mathbb{G},(\widetilde{W},L))$.
\ethe

\brem It is proved in \cite{JS} that, for any probability measure $\mathbb{P}$ on $\mathcal{F}_\infty$, for any positive $(\mathbb{P},\mathbb{F})$ local martingale $N$, for any continuous $\mathbb{F}$-adapted increasing process $\Lambda$ such that $\Lambda_0=0, N_0=1$ and $\forall t>0,N_{t}e^{-\Lambda_t}< 1, N_{t-}e^{-\Lambda_t}< 1$, there exists a probability measure $\mathbb{Q}$ on the product space, which is absolutely continuous with respect to the Cox measure, such that $\mathbb{Q}|_{\mathcal{F}_\infty}=\mathbb{P}|_{\mathcal{F}_\infty}$ and $\mathbb{Q}[t<\tau|\mathcal{F}_\infty]=N_te^{-\Lambda_t}, t\geq 0$.  
\erem

\section{An evolution model}\label{JSmodel}

In this section, we consider a model developed in \cite{JS2} (the $\natural$-model). The basic setting is a filtered probability space $(\Omega,\mathcal{A},\mathbb{F},\mathbb{P})$, where
$\mathbb{F}=(\mathcal{F}_t)_{t\geq 0}$ is a
filtration satisfying the usual conditions. We consider an $\mathbb{F}$-adapted continuous
increasing process $\Lambda$ and a non negative $(\mathbb{P},\mathbb{F})$ local martingale $N$ such that $\Lambda_0=0,N_0=1$ and $0\leq N_t e^{-\Lambda_t}\leq 1$ for all $0\leq t<\infty$. We introduce the product measurable space $([0,\infty]\times\Omega,
\mathcal{B}[0,\infty]\otimes\mathcal{F}_\infty)$ with its canonical projection maps $\pi$ and $\tau$ : $\pi(s,\omega)=\omega$ and $\tau(s,\omega)=s$. Using the projection map $\pi$ we pull back the probability structure $(\mathbb{P},\mathbb{F})$ onto the product space $[0,\infty]\times\Omega$ with the filtration $\hat{\mathbb{F}}=\pi^{-1}(\mathbb{F})$ and with the probability measure on $\pi^{-1}(\mathcal{F}_\infty)$ defined by $\hat{\mathbb{P}}(\pi^{-1}(\mathtt{A}))=\mathbb{P}(\mathtt{A})$ for $\mathtt{A}\in\mathcal{F}_\infty$. The probability structure $([0,\infty]\times\Omega,\hat{\mathbb{F}},\hat{\mathbb{P}})$ is isomorphic to that of $(\Omega,\mathbb{F},\mathbb{P})$. We will henceforth simply denote $(\hat{\mathbb{P}},\hat{\mathbb{F}})$ by $(\mathbb{P},\mathbb{F})$ and identify the $\mathcal{F}_\infty$-measurable random variables $\xi$ on $\Omega$ with $\xi\circ\pi$ on the product space. 

We consider the following problem:

\noindent\textbf{Problem $\mathcal{P}^*$.} Construct on the
product space $([0,\infty]\times\Omega,
\mathcal{B}[0,\infty]\otimes\mathcal{F}_\infty)$ a probability measure
$\Q$ such that 
\begin{itemize}
\item
(restriction condition) $\mathbb{Q}|_{\mathcal{F}_\infty}=\mathbb{P}|_{\mathcal{F}_\infty}$ and
\item
(projection condition) $\mathbb{Q}[\tau>t|\mathcal{F}_t]=N_t
e^{-\Lambda_t}$ for all $0\leq t<\infty$.
\end{itemize} 
(Recall that we identify $\hat{\mathbb{F}}$ as $\mathbb{F}$ and $\hat{\mathbb{P}}$ as $\mathbb{P}$.)

Suppose \textbf{Hy}(C), i.e. all $(\mathbb{P},\mathbb{F})$ local martingales are continuous. Suppose $Z_t< 1$ for any $0<t<\infty$, where $Z=Ne^{-\Lambda}$. Under these conditions, \cite{JS2} proves that there exist infinitely many solutions to the problem $\mathcal{P}^*$. In particular, for any $(\mathbb{P},\mathbb{F})$ local martingale $Y$, for any bounded differentiable function $f$ with bounded continuous derivative and $f(0)=0$, there exists $\mathbb{Q}^\natural$ a solution of the problem $\mathcal{P}^*$ on the product space such that, for any $u\in\mathbb{R}_+^*$, the martingale $M^u_t=\mathbb{Q}^\natural[\tau\leq u|\mathcal{F}_t], t\geq u,$ satisfies the following evolution equation($\natural$) $$
(\natural_u)  \left\{\dcb
dX_t=X_t\left(-\frac{e^{-\Lambda_t}}{1-Z_t}dN_t+f(X_t - (1-Z_t))dY_t\right),\ u\leq t<\infty,\\
X_u=1-Z_u.
\dce
\right.
$$
Consider the progressively enlarged filtration $\mathbb{G}$ on this product space with the random time $\tau$.

\bethe \label{decompositionFormula}
Suppose the same assumptions as above. Suppose in addition : 

\noindent\textbf{Hy}(Mc) : For each $0< t< \infty$, the map $u\rightarrow M^u_t$ is continuous on $(0,t]$.

\noindent Then, for any $(\mathbb{P},\mathbb{F})$ local martingale $X$, the process 
$$
\dcb
\pv(X)_t:=\int_{0}^{t}\ind_{\{s\leq \tau\}}\frac{e^{-\Lambda_s}}{Z_s}d\cro{N,X}_s
-\int_{0}^{t}\ind_{\{\tau< s\}}\frac{e^{-\Lambda_s}}{1-Z_s}d\cro{N,X}_s\\
\\
\hspace{1.5cm}+\int_{0}^{t}\ind_{\{\tau< s\}}( f(M^{\tau}_s-(1-Z_s))+M^\tau_s f'(M^{\tau}_s-(1-Z_s))) d\cro{Y,X}_s,\ \ 0\leq t<\infty,
\dce
$$is a well-defined $\mathbb{G}$-predictable process with finite variation, and the process
$
\widetilde{X}=X-\pv(X)
$ is a $(\mathbb{Q}^\natural,\mathbb{G})$ local martingale. 
\ethe

We now study the \mrt\ property of the model defined by the evolution equation$(\natural)$. Let $W$ be a $d$-dimensional $\mathbb{F}$-adapted càdlàg process. We assume the following set of assumptions :

\bassumption\label{assumptionDieze}
\ebe
\item[(i)]
The above two parameters $Y$ and $f$ are given. 
\item[(ii)]
Assume \textbf{Hy}(C), \textbf{Hy}(Mc).
\item[(iii)]
Assume \mrt$(\mathbb{P},\mathbb{F},W)$. 
\item[(iv)]
Assume $0<Z_t<1$ for $0<t<\infty.$
\dbe
\eassumption

For any $\mathbb{F}$ stopping time $T$, under assumption \textbf{Hy}(Mc), we have$$
\mathbb{Q}^\natural[\tau=T, \tau\leq t|\mathcal{F}_t]=\int_0^t\ind_{\{u=T\}}d_uM^u_t=0,\ \forall 0<t<\infty.
$$
This yields that $\tau$ avoids the $\mathbb{F}$ stopping times on $(0,\infty)$.
Let $0<a<\infty,n\in\mathbb{N}^\ast$ and let $$
\dcb
T_{a,n}&=\inf\{v\geq a:\hspace{3pt} & \int_a^v\frac{e^{-2\Lambda_w}}{Z_w^2}d\cro{N}_w>n,\\
&&\mbox{ or }\int_a^v\frac{e^{-2\Lambda_w}}{(1-Z_w)^2}d\cro{N}_w>n,\\
&&\mbox{ or }\cro{Y}_v-\cro{Y}_a>n,\\
&&\mbox{ or }\cro{W}_v-\cro{W}_a>n,\\
&&\mbox{ or }v>a+n \hspace{2cm}  \}.
\dce
$$
The random times $T_{a,n}$ are $\mathbb{F}$-stopping times. Since $0<Z<1$ on $(0,\infty)$ and $N,Y,W$ are continuous, $\lim_{n\rightarrow\infty}T_{a,n}=\infty$. We have $(0,\infty)=\cup_{a\in\mathtt{Q},n\in\mathbb{N}^*}(a,T_{a,n})$. 

Let us show that there exists an $s\!\mathcal{H}$-measure on the intervals $(a,T_{a,n}]$. We introduce $$
\gamma_s=\frac{e^{-\Lambda_s}}{Z_s},\
\alpha_s=-\frac{e^{-\Lambda_s}}{1-Z_s},\
\beta_s=f(M^{\tau}_s-(1-Z_s))+M^\tau_s f'(M^{\tau}_s-(1-Z_s))
$$
and the exponential martingale: $$
\eta=\mathcal{E}\left((-\gamma\ind_{[0,\tau]}-\alpha\ind_{(\tau,\infty)})\ind_{(a,T_n]}\centerdot \widetilde{N}
+(-\beta)\ind_{(\tau,\infty)}\ind_{(a,T_n]}\centerdot \widetilde{Y}\right).
$$
By \textbf{Hy}(C), $\pv(N)$ is continuous and $\cro{N}=[N]=\cro{\widetilde{N}}$ indifferently in the filtration $\mathbb{F}$ or in the filtration $\mathbb{G}$. The same property holds for the bracket of $Y$. We check then that Novikov's condition is satisfied by $\eta$ so that $\mathbb{Q}^\natural[\eta]=1$. Let $
\mathbb{Q}^{a,n}=\eta\cdot \mathbb{Q}^\natural.
$
By Girsanov's theorem, the process $\widetilde{X}^{[\eta]}_{ t}-\widetilde{X}^{[\eta]}_{a},\ 0\leq t<\infty$, is a $(\mathbb{Q}^{a,n},\mathbb{G})$ local martingale. Note$$
\pv(X)=\gamma\ind_{[0,\tau]}\centerdot \cro{N,X}
+\alpha\ind_{(\tau,\infty)}\centerdot \cro{N,X}
+\beta\ind_{(\tau,\infty)}\centerdot \cro{Y,X}.
$$ 
Because of \textbf{Hy}(C), we can write $\cro{\widetilde{X},\widetilde{N}}=\cro{X,N}$ and $\cro{\widetilde{X},\widetilde{Y}}=\cro{X,Y}$, and therefore, by a direct computation (cf. subsection \ref{honest}), we get
$$
\widetilde{X}^{[\eta]}_{t}-\widetilde{X}^{[\eta]}_{a}
=X_{t}-X_{a},\ \mbox{ if $a\leq t\leq T_{a,n}$}.
$$ 
This shows that $X^{T_{a,n}}_{a \vee t}-X_{a}$ is a $(\mathbb{Q}^{a,n},\mathbb{G})$ local martingale. Since $X_{T_{a,n}}-X_{a}\in\mathcal{F}_{T_{a,n}}\subset \mathcal{G}^*_{T_{a,n}}$, $X^{(a,T_{a,n}]}$ is also a $(\mathbb{Q}^{a,n},\mathbb{G}^{(a,T_{a,n}]})$ local martingale (cf. Proposition \ref{appendix} (5)). The measure $\mathbb{Q}^{a,n}$ is an $s\!\mathcal{H}$-measure on $(a, T_{a,n}]$. 

The $s\!\mathcal{H}$-measure condition covering $(0,\infty)$ is satisfied. Applying Theorem \ref{G=G}, Lemma \ref{reducing}, Theorem \ref{put-together}, we obtain

\bethe\label{diez}
Under Assumption \ref{assumptionDieze}, the property \mrt$(\mathbb{Q}^\natural,\mathbb{G},(\widetilde{W},L))$ holds.
\ethe

\

\

\appendix

\section{Study of the filtration $\mathbb{G}^{(S,T]}$}

In this appendix, we study the filtrations $\mathbb{G}^{(S,T]}$ introduced in Section \ref{frags} and prove the results stated in Proposition \ref{appendix}. This study is independent of the main text of this article. We consider a probability space $(\Omega,\mathcal{A},\mathbb{Q})$, with a right-continuous filtration $\ff=(\mathcal{F}_t)_{t\geq 0}$ of sub-$\sigma$-algebras in $\mathcal{A}$. We consider a random time $\tau$ in $\mathcal{A}$ and its associated progressively enlarged filtration $\mathbb{G}=(\mathcal{G}_t)_{t\geq 0}$ where $\G_t=\cap_{s>t}(\F_{s}\vee \sigma (\tau \wedge s))$ for $t\in\mathbb{R}_+$. Unlike the assumptions in Subsection \ref{framework}, we do not assume here that $\mathcal{G}_0$ contains all the $(\mathbb{Q},\mathcal{G}_\infty)$ negligible sets.  

We recall two useful results:

\bl\label{partTribu}
Let $E$ be a space. Let $\mathcal{C}$ be a non empty family of sets in $E$. Let $A\subset E$. Then, $A\cap \sigma(\mathcal{C})=A\cap\sigma(A\cap \mathcal{C})$. 
\el

\bl\label{swing}
(cf. \cite[Lemme (4.4)]{J}) Let $T$ be a $\mathbb{G}$-stopping time. We have $$
\dcb
\mathcal{G}_{T-}&=& \{T\leq\tau,T<\infty\}\cap\mathcal{F}_{T-}+\{\tau< T<\infty\}\cap(\sigma(\tau)\vee\mathcal{F}_{T-})+\{T=\infty\}\cap(\sigma(\tau)\vee\mathcal{F}_{\infty})\\

&=& \{T<\tau\}\cap\mathcal{F}_{T-}+\{\tau\leq T<\infty\}\cap(\sigma(\tau)\vee\mathcal{F}_{T-})+\{T=\infty\}\cap(\sigma(\tau)\vee\mathcal{F}_{\infty}).
\dce
$$
\el

From now on till Lemma \ref{gSTKK}, we consider two $\mathbb{G}$-stopping times $S,T$ such that $S\leq T$. 

\bl\label{gSTcharact}
We have the relations : $\mathcal{G}^{(S,T]}_\infty=\mathcal{G}^*_T\subset\sigma(\tau)\vee\mathcal{F}_T$ and $\mathcal{G}_{T-}\subset \mathcal{G}^\ast_{T}\subset \mathcal{G}_{T}$. For $t\geq 0$, we have $\mathcal{G}^{(S,T]}_t\subset \mathcal{G}_{(S\vee t)\wedge T}$. A random variable $\xi$ is $\mathcal{G}^{(S,T]}_t$-measurable if and only if there exist two random variables $\xi_1\in\mathcal{G}^*_T$ and $\xi_2\in\mathcal{G}_{S\vee t}$ such that $$
\xi = \xi_1\ind_{\{T\leq S\vee t\}}+\xi_2\ind_{\{S\vee t<T\}}
$$or if and only if $\xi\ind_{\{T\leq S\vee t\}}\in\mathcal{G}^*_T$ and $\xi\ind_{\{S\vee t<T\}}\in\mathcal{G}_{S\vee t}$.
\el

\proof The results of the lemma are direct consequences of the definition. Let us prove, for example, the relation $\mathcal{G}^*_T\subset\sigma(\tau)\vee\mathcal{F}_T$. We note first that $\ind_{[0,t]}(T)\in\mathcal{F}_{T-}\subset \mathcal{F}_{T}$, and consequently, $$
\{T<\tau\} = \cup_{t\in\mathtt{Q}_+}\{T\leq t\}\cap\{t<\tau\}\in\sigma(\tau)\vee\mathcal{F}_T.
$$
In the same way, we can prove that $\{\tau\leq T\}$ and $\{ T=\infty\}$ belong to $\sigma(\tau)\vee\mathcal{F}_T$. Therefore,$$
\dcb
\mathcal{G}^\ast_{T}
&\subset&\{T<\tau\}\cap(\sigma(\tau)\vee\mathcal{F}_{T})+\{\tau\leq T<\infty\}\cap(\sigma(\tau)\vee\mathcal{F}_{T})+\{T=\infty\}\cap(\sigma(\tau)\vee\mathcal{F}_{T})\\

&=&\sigma(\tau)\vee\mathcal{F}_{T},
\dce
$$
where we have used the relation $\{T=\infty\}\cap\mathcal{F}_{T}=\{T=\infty\}\cap\mathcal{F}_{\infty}$. The lemma is proved.
\ok

\bl(Proposition \ref{appendix} (1))
$\mathbb{G}^{(S,T]}$ is a right-continuous filtration.

\el

\proof Let us show first that $\mathbb{G}^{(S,T]}$ is a filtration, i.e., $\mathcal{G}^{(S,T]}_s\subset \mathcal{G}^{(S,T]}_t$ for any $0<s<t<\infty$. By definition,
$$
\dcb
&&\mathcal{G}^{(S,T]}_t \\

&=&\{\{T\leq S\vee t\}\cap(\{T\leq S\vee s\}\cap A + \{S\vee s<T\}\cap A) + \{S\vee t<T\}\cap B : A\in\mathcal{G}^\ast_T,\ B\in\mathcal{G}_{S\vee t}\}\\

&=&\{\{T\leq S\vee s\}\cap A + \{S\vee s<T\}\cap (\{T\leq S\vee t\}\cap A' + \{S\vee t<T\}\cap B) : A,A'\in\mathcal{G}^\ast_T,\ B\in\mathcal{G}_{S\vee t}\}.\\

\dce
$$
We study $\mathcal{G}^{(S,T]}$ separately in two cases : $\{t\leq S\}$ or $\{S<t\}$. Firstly, we have$$
\dcb
&&\{t\leq S\}\cap \mathcal{G}^{(S,T]}_t\\

&=&\{\{t\leq S\}\cap\{T\leq S\}\cap A + \{t\leq S\}\cap\{S<T\}\cap B : A\in\mathcal{G}^\ast_T,\ B\in\mathcal{G}_{S}\}\\

&=&\{\{t\leq S\}\cap\{T\leq S\vee s\}\cap A + \{t\leq S\}\cap\{S\vee s<T\}\cap B : A\in\mathcal{G}^\ast_T,\ B\in\mathcal{G}_{S\vee s}\}\\

&=&\{t\leq S\}\cap \mathcal{G}^{(S,T]}_{s}.
\dce
$$
Next, we note that, on the set $\{S<t\}$, we have $S\vee s <t$. Therefore
$$
\dcb
&&\{S<t\}\cap \mathcal{G}^{(S,T]}_t\\

&=&\{\{S<t\}\cap\{T\leq S\vee s\}\cap A + \{S<t\}\cap\{S\vee s<T\}\cap (\{T\leq t\}\cap A' + \{t<T\}\cap B) : \\
&&\hspace{12cm}A,A'\in\mathcal{G}^\ast_T,\ B\in\mathcal{G}_{S\vee t}\}\\

&\supset&\{\{S<t\}\cap\{T\leq S\vee s\}\cap A + \{S<t\}\cap\{S\vee s<T\}\cap (\{T\leq t\}\cap A' + \{t<T\}\cap B) : \\
&&\hspace{11cm}A\in\mathcal{G}^*_{T},\ A'\in\mathcal{G}_{T-},\ B\in\mathcal{G}_{t-}\}\\

&&\mbox{ \hspace{1cm} because, from Lemma \ref{gSTcharact}, $\mathcal{G}_{T-} \subset \mathcal{G}^*_{T}$,}\\
\\
&=&\{\{S<t\}\cap\{T\leq S\vee s\}\cap A + \{S<t\}\cap\{S\vee s<T\}\cap C : A\in\mathcal{G}^*_{T},\ C\in\mathcal{G}_{(T\wedge t)-}\}\\

&&\mbox{\hspace{1cm}  see \cite[Corollary 3.5]{Yan}}\\

&\supset&\{\{S<t\}\cap\{T\leq S\vee s\}\cap A + \{S<t\}\cap\{S\vee s<T\wedge t\}\cap C : A\in\mathcal{G}^*_{T},\ C\in\mathcal{G}_{S\vee s}\}\\

&=&\{S<t\}\cap\mathcal{G}^{(S,T]}_s. \\
\dce
$$
This shows that $\mathbb{G}^{(S,T]}$ is a filtration. 

Now we show the right-continuity, i.e., $\mathcal{G}^{(S,T]}_s=\cap_{t>s}\mathcal{G}^{(S,T]}_t$ for any $s\in\mathbb{R}_+$. We begin with the following observation. Let $D$ be a set in $\Omega$. We have the equality \begin{equation}\label{DG1}
D\cap(\cap_{t\in\mathtt{Q},t>s} \mathcal{G}^{(S,T]}_t)
=\cap_{t\in\mathtt{Q},t>s}(D\cap \mathcal{G}^{(S,T]}_t).
\end{equation} 
Clearly the left-hand side family is contained in the right one. Let $B$ be an element in the right-hand side family. For any $t\in\mathtt{Q}$ with $t>s$, there exists $B_t\in\mathcal{G}^{(S,T]}_t$ such that $B=D\cap B_t$. It yields that, for $\beta>0$, $B=\cup_{t\in\mathtt{Q},\beta>t>s}(D\cap B_t)=D\cap(\cup_{t\in\mathtt{Q},\beta>t>s} B_t)$ and $$
B=\cap_{\beta>s}\cup_{t\in\mathtt{Q},\beta>t>s}(D\cap B_t)=D\cap (\cap_{\beta>s}\cup_{t\in\mathtt{Q},\beta>t>s}B_t).
$$ 
By monotony $$
\cap_{\beta>s}\cup_{t\in\mathtt{Q},\beta>t>s}B_t
=\cap_{s+\epsilon>\beta>s}\cup_{t\in\mathtt{Q},\beta>t>s}B_t
\in\mathcal{G}^{(S,T]}_{s+\epsilon}
$$
for any $\epsilon>0$. This shows that $B$ is an element in the left-hand side family. The formula (\ref{DG1}) is proved. In the same way we prove also the equality $$
D\cap(\cap_{t\in\mathtt{Q},t>s} \mathcal{G}_{S\vee t})=\cap_{t\in\mathtt{Q},t>s} (D\cap\mathcal{G}_{S\vee t}).
$$
Using the above observation, we obtain, on the one hand, $$
\dcb
\{T\leq S\vee s\}\cap(\cap_{t>s} \mathcal{G}^{(S,T]}_t)
&=& \cap_{t\in\mathtt{Q},t>s}(\{T\leq S\vee s\}\cap \mathcal{G}^{(S,T]}_t) \\

&=&\{T\leq S\vee s\}\cap \mathcal{G}^*_T\\
&=& \{T\leq S\vee s\}\cap \mathcal{G}^{(S,T]}_s,
\dce
$$
and on the other hand, for $\beta>s$,$$
\dcb
\{S\vee \beta<T\}\cap(\cap_{t>s} \mathcal{G}^{(S,T]}_t)

&=&\cap_{t\in\mathtt{Q},\beta>t>s} (\{S\vee \beta<T\}\cap\mathcal{G}^{(S,T]}_t)\\

&=&\{S\vee \beta<T\}\cap(\cap_{t\in\mathtt{Q},\beta>t>s} \mathcal{G}_{S\vee t})\\

&=&\{S\vee \beta<T\}\cap \mathcal{G}^{(S,T]}_s.\\

\dce
$$
Consider then a set $B\in\{S\vee s<T\}\cap(\cap_{t>s} \mathcal{G}^{(S,T]}_t)$. For any $\beta>s$, $$
B\cap\{S\vee \beta<T\}\in\{S\vee \beta<T\}\cap(\cap_{t>s} \mathcal{G}^{(S,T]}_t)
=\{S\vee \beta<T\}\cap \mathcal{G}^{(S,T]}_s.
$$
Let $B_\beta\in\mathcal{G}^{(S,T]}_s$ such that $B\cap\{S\vee \beta<T\}=B_\beta\cap\{S\vee \beta<T\}$. We have $B_\beta\cap\{S\vee \beta<T\}=B_{\beta'}\cap\{S\vee \beta<T\}$ for any $\beta>\beta'>s$. Let $B_*=\cap_{\epsilon>0}\cup_{\beta'\in\mathtt{Q},s<\beta'<s+\epsilon}B_{\beta'}\in\mathcal{G}^{(S,T]}_s$. We have$$
\dcb
B_*\cap\{S\vee \beta<T\}
&=&\{S\vee \beta<T\}\cap(\cap_{\epsilon>0}\cup_{s<\beta'<s+\epsilon}B_{\beta'})\\
&=&\{S\vee \beta<T\}\cap(\cap_{\beta-s>\epsilon>0}\cup_{s<\beta'<s+\epsilon}B_{\beta'})\\
&=&\cap_{\beta-s>\epsilon>0}\cup_{s<\beta'<s+\epsilon}(\{S\vee \beta<T\}\cap B_{\beta'})\\
&=&\cap_{\beta-s>\epsilon>0}\cup_{s<\beta'<s+\epsilon}(\{S\vee \beta<T\}\cap B)\\
&=&B\cap \{S\vee \beta<T\}.
\dce
$$
Taking the union on $\beta>s$ we obtain finally$$
B=B\cap\{S\vee s<T\}=B_*\cap\{S\vee s<T\}.
$$
This being true for any $B\in\{S\vee s<T\}\cap(\cap_{t>s} \mathcal{G}^{(S,T]}_t)$,
 we obtain 
$$
\{S\vee s<T\}\cap(\cap_{t>s} \mathcal{G}^{(S,T]}_t)\subset\{S\vee s<T\}\cap \mathcal{G}^{(S,T]}_s.
$$ 
Actually we have an equality instead of the inclusion, because the inverse of the above relation is an evidence. Since $\{S\vee s<T\}\in \mathcal{G}^{(S,T]}_s$, we can put together the two equalities: 
$$
\dcb
\{T\leq S\vee s\}\cap(\cap_{t>s} \mathcal{G}^{(S,T]}_t)&=&\{T\leq S\vee s\}\cap \mathcal{G}^{(S,T]}_s\\
\{S\vee s<T\}\cap(\cap_{t>s} \mathcal{G}^{(S,T]}_t)&=&\{S\vee s<T\}\cap \mathcal{G}^{(S,T]}_s.
\dce
$$ 
and we conclude $\mathcal{G}^{(S,T]}_s=\cap_{t>s}\mathcal{G}^{(S,T]}_t$, i.e., the right continuity of $\mathbb{G}^{(S,T]}$. \ok

\bl\label{gSTstoppingtimes}
$S$ and $T$ are $\mathbb{G}^{(S,T]}$-stopping times. More generally, any $\mathbb{G}$-stopping time $R$ between $S$ and $T$ : $S\leq R\leq T$, is a $\mathbb{G}^{(S,T]}$-stopping time. Conversely, for any $\mathbb{G}^{(S,T]}$-stopping time $R'$, $S\vee R'$ is a $\mathbb{G}$-stopping time.
\el

This lemma can be checked straightforwardly.

\bl\label{gSTgST}(Proposition \ref{appendix} (2))
For a $\mathbb{G}^{(S,T]}$-stopping time $R$ such that $S\leq R\leq T$, we have $$
\mathcal{G}^{(S,T]}_R=\{R=T\}\cap\mathcal{G}^*_{T}+\{R<T\}\cap\mathcal{G}_{R}
=\{R=T\}\cap\mathcal{G}^*_{R}+\{R<T\}\cap\mathcal{G}_{R}.
$$
In particular, $\mathcal{G}^{(S,T]}_S=\mathcal{G}^{(S,T]}_0$ and $\mathcal{G}^{(S,T]}_T=\mathcal{G}^*_{T}$. 
\el

\proof Consider first $\mathcal{G}^{(S,T]}_T$. We have, on the one hand, $\mathcal{G}^{(S,T]}_T\subset \mathcal{G}^{(S,T]}_\infty= \mathcal{G}^*_T$. On the other hand, for any $B\in \mathcal{G}^*_{T}$, for $t\geq 0$, we note that $\{T\leq t\}\in\mathcal{G}_{T-}\subset \mathcal{G}^*_{T}$ so that $B\cap \{T\leq t\}\in \mathcal{G}^*_{T}$. Therefore,
$$
\dcb
B\cap \{T\leq t\}
&=&B\cap \{T\leq t\}\cap\{T\leq S\vee t\} + B\cap \{T\leq t\}\cap\{S\vee t < T\}\\

&=&B\cap \{T\leq t\}\cap\{T\leq S\vee t\}

\in\{T\leq S\vee t\}\cap \mathcal{G}^*_{T}

\in \mathcal{G}^{(S,T]}_{t}.
\dce
$$
This proves that $B\in \mathcal{G}^{(S,T]}_T$, i.e., $\mathcal{G}^{(S,T]}_T=\mathcal{G}^*_{T}$.

Consider now $\mathcal{G}^{(S,T]}_R$. We note that $R$ is also a $\mathbb{G}$-stopping time, as it is stated in the preceding Lemma \ref{gSTstoppingtimes}. Let $A\in \mathcal{G}^{(S,T]}_R$. Then, for any $t\geq 0$, since $\{R<T\}$ is an element in $\mathcal{G}^{(S,T]}_R$, we have $$
\dcb
&&A\cap \{R<T\}\cap \{R\leq t\}
\in \{R\leq t\}\cap\mathcal{G}^{(S,T]}_t\\
&\subset&  \{R\leq t\}\cap\mathcal{G}_{(S\vee t)\wedge T}\ \mbox{ according to Lemma \ref{gSTcharact}}\\

&\subset&  \{R\leq t\}\cap\{S\leq t\}\cap\mathcal{G}_{S\vee t}

\subset\mathcal{G}_{t}.
\dce
$$ 
This computation shows that $A\cap \{R<T\}\in\mathcal{G}_R$. Conversely let $B\in\mathcal{G}_R$. For $t\geq 0$, we have $$
B\cap \{R<T\}\cap \{R\leq t\}\in\mathcal{G}_{T-} \mbox{ and }
B\cap \{R<T\}\cap \{R\leq t\}\in\mathcal{G}_{t}. 
$$ 
Therefore,
$$
\dcb
&&B\cap \{R<T\}\cap \{R\leq t\}\\

&=&\{T\leq S\vee t\}\cap B\cap \{R<T\}\cap \{R\leq t\}+\{S\vee t<T\}\cap B\cap \{R<T\}\cap \{R\leq t\}\\

&\in&\{T\leq S\vee t\}\cap \mathcal{G}_{T-}+\{S\vee t<T\} \cap \{R\leq t\}\cap \mathcal{G}_{t}\\

&\subset&\{T\leq S\vee t\}\cap \mathcal{G}^*_{T}+\{S\vee t<T\} \cap \mathcal{G}_{S\vee t}\\

&=&  \mathcal{G}^{(S,T]}_{t}.
\dce
$$
This proves that $B\cap \{R<T\}\in \mathcal{G}^{(S,T]}_R$. We have just established $$
\{R<T\}\cap \mathcal{G}^{(S,T]}_R=\{R<T\}\cap\mathcal{G}_{R}.
$$We can also write $$
\{R=T\}\cap \mathcal{G}^{(S,T]}_R=\{R=T\}\cap \mathcal{G}^{(S,T]}_T=\{R=T\}\cap\mathcal{G}^*_{T}.
$$ 
(For the first equality, see \cite[Corollary 3.5 statement 4)]{Yan}). The lemma is proved. \ok

\bl\label{gSTandG}(Proposition \ref{appendix} (3))
For any $\mathbb{G}$-adapted process $X$ such that $X_T\in\mathcal{G}^*_T$, $X^{(S,T]}$ defines a $\mathbb{G}^{(S,T]}$-adapted process. Conversely, for any $\mathbb{G}^{(S,T]}$-adapted process $X'$, $X'^{(S,T]}$ defines a $\mathbb{G}$-adapted process.
\el

This lemma can be checked using Lemma \ref{gSTcharact}.

\bl\label{gSTmartingales}(Proposition \ref{appendix} (5))
For any $(\mathbb{Q},\mathbb{G})$ local martingale $X$ such that $X_T\in\mathcal{G}^*_T$, $X^{(S,T]}$ defines a $(\mathbb{Q},\mathbb{G}^{(S,T]})$ local martingale. Conversely, for any $(\mathbb{Q},\mathbb{G}^{(S,T]})$ local martingale $X'$, $X'^{(S,T]}$ defines a $(\mathbb{Q},\mathbb{G})$ local martingale.
\el

\proof Let $R$ be a $X$-reducing stopping time, i.e., a $\mathbb{G}$-stopping time such that $X^R$ is a uniformly $\mathbb{Q}$-integrable $(\mathbb{Q},\mathbb{G})$ martingale. Note that, for $t\geq 0$, $
(X^{(S,T]})^R_t =(X^R)^{(S,T]}_t.
$
As $\ind_{\{R<T\}}X_R,\ind_{\{R\geq T\}}\in\mathcal{G}_{T-}$, we have $$
X^R_T=\ind_{\{R<T\}}X_R+\ind_{\{R\geq T\}}X_T \in\mathcal{G}^*_{T}.
$$
This relation together with the preceding Lemma \ref{gSTandG} shows that $(X^{(S,T]})^R$ is $\mathbb{G}^{(S,T]}$-adapted. (We note that we could not say this directly, because $R$ is not a $\mathbb{G}^{(S,T]}$-stopping time.)

Let $0\leq s \leq t \leq \infty$ and $A\in\mathcal{G}^{(S,T]}_s$. Since $\mathcal{G}^{(S,T]}_s\subset \mathcal{G}_{(S\vee s)\wedge T}\subset \mathcal{G}_{(S\vee s)}$, we have$$
\dcb
\mathbb{E}[\ind_A (X^R)^T_{S\vee t}] = \mathbb{E}[\ind_A X^{R\wedge T}_{S\vee t}]
=\mathbb{E}[\ind_AX^{R\wedge T}_{S\vee s}]=\mathbb{E}[\ind_A (X^R)^T_{S\vee s}].\\
\dce
$$
This computation implies $$
\mathbb{E}[\ind_A (X^R)^{(S,T]}_{S\vee t}]
=\mathbb{E}[\ind_A (X^R)^{(S,T]}_{S\vee s}],
$$ 
i.e. $(X^R)^{(S,T]}=(X^{(S,T]})^R$ is a $(\mathbb{Q},\mathbb{G}^{(S,T]})$ uniformly integrable martingale. We note that $X^{(S,T]}_t=0$ if $t\leq S$ and $X^{(S,T]}_t=X_T-X_S$ if $t\geq T$. We can write$$
(X^{(S,T]})^R  = (X^{(S,T]})^{R'},
$$
where $$
\dcb
R'
&=&\ind_{\{(S\vee R)\wedge T< T\}}(S\vee R)\wedge T+\ind_{\{(S\vee R)\wedge T\geq  T\}}\cdot \infty.
\dce
$$
According to Lemma \ref{gSTstoppingtimes}, $(S\vee R)\wedge T$, as well as $T$, is a $\mathbb{G}^{(S,T]}$-stopping time. The set $\{(S\vee R)\wedge T< T\}$ is in $\mathcal{G}^{(S,T]}_{(S\vee R)\wedge T}$. Hence, $R'$, as a restriction of $(S\vee R)\wedge T$, is also a $\mathbb{G}^{(S,T]}$-stopping time. Moreover, if $R$ tends to infinity, $R'$ does so too. We can now state that $X^{(S,T]}$ is a $(\mathbb{Q},\mathbb{G}^{(S,T]})$ local martingale. The first part of the lemma is proved. 

Consider a $(\mathbb{Q},\mathbb{G}^{(S,T]})$ local martingale $X'$. Let $U'$ be a $X'$-reducing stopping time. We know that $X'^{(S,T]}$ is a $\mathbb{G}$-adapted process, and $S\vee U'$ is a $\mathbb{G}$-stopping time. We use the identity $$
\dcb
(X'^{(S,T]})^{S\vee U'}_t
&=&(X'^{S\vee U'})^{(S,T]}_t=(X'^{S\vee U'})^T_{S\vee t}-X'^{S\vee U'}_S=X'^{(S\vee U')\wedge T}_{S\vee t}-X'_S
\dce
$$ 
for $t\geq 0$. Set $R=(S\vee U')\wedge T$ which is a stopping time with respect to both $\mathbb{G}$ and $\mathbb{G}^{(S,T]}$. Let $0\leq s\leq t\leq \infty$ and $B\in\mathcal{G}_s$. We write$$
\dcb
\mathbb{E}[\ind_B (X'^{(S,T]})^{S\vee U'}_t] 
&=& \mathbb{E}[\ind_B\ind_{\{T\leq S\vee s\}} (X'^{R}_{S\vee t}-X'_S)] + \mathbb{E}[\ind_B\ind_{\{S\vee s<T\}} (X'^{R}_{S\vee t}-X'_S)].\\
\dce
$$
We note that $B\cap\{S\vee s<T\}\in\mathcal{G}_{S\vee s}$ which yields $B\cap\{S\vee s<T\}\in\mathcal{G}^{(S,T]}_s$. Hence, $$
\dcb
&&\mathbb{E}[\ind_B\ind_{\{S\vee s<T\}} (X'^{R}_{S\vee t}-X'_S)]
=\mathbb{E}[\ind_B\ind_{\{S\vee s<T\}} \ind_{\{S<U'\}} (X'^{U'\wedge T}_{S\vee t}-X'_S)]\\
&=&\mathbb{E}[\ind_B\ind_{\{S\vee s<T\}} \ind_{\{S<U'\}} (X'^{U'\wedge T}_{S\vee s}-X'_S)]\
\mbox{ because $\{S<U'\}\in\mathcal{G}^{(S,T]}_{S\vee s}$}\\
&=&\mathbb{E}[\ind_B\ind_{\{S\vee s<T\}} (X'^{R}_{S\vee s}-X'_S)].
\dce
$$
We note also $$
\ind_B\ind_{\{T\leq S\vee s\}} (X'^{R}_{S\vee t}-X'_S)
=\ind_B\ind_{\{T\leq S\vee s\}} (X'_{R}-X'_S)
=\ind_B\ind_{\{T\leq S\vee s\}} (X'^R_{S\vee s}-X'_S).
$$
Putting these relations together, we write $$
\dcb
\mathbb{E}[\ind_B (X'^{(S,T]})^{S\vee U'}_t] 
&=& \mathbb{E}[\ind_B\ind_{\{T\leq S\vee s\}} (X'^{R}_{S\vee s}-X'_S)] + \mathbb{E}[\ind_B\ind_{\{S\vee s<T\}} (X'^{R}_{S\vee s}-X'_S)]\\
&=&\mathbb{E}[\ind_B (X'^{R}_{S\vee s}-X'_S)]
=\mathbb{E}[\ind_B (X'^{(S,T]})^{S\vee U'}_t].
\dce
$$This proves the second part of the lemma. 
\ok

\bl\label{dSTpredictable}(Proposition \ref{appendix} (4))
For any $\mathbb{G}$-predictable process $K$, the processes $K^{(S,T]}$ and $\ind_{(S,T]}K$ are $\mathbb{G}^{(S,T]}$-predictable. Conversely, for any $\mathbb{G}^{(S,T]}$-predictable process $K'$, $K'^{(S,T]}$ and $\ind_{(S,T]}K'$ are $\mathbb{G}$-predictable processes.
\el

\proof If $K$ or $K'$ are left continuous, with Lemma \ref{gSTandG}, the result is clear. To pass to general $K$ or $K'$, we use the monotone class theorem. \ok

\bl\label{gSTVV}
For any increasing process $A$, $(\ind_{(S,T]})\stocint A$ is $(\mathbb{Q},\mathbb{G})$-locally integrable if and only if $(\ind_{(S,T]})\stocint A$ is $(\mathbb{Q},\mathbb{G}^{(S,T]})$-locally integrable. In the case of a locally integrable $A$, let $A''$ and $A'$ be respectively the $(\mathbb{Q},\mathbb{G})$ and the $(\mathbb{Q},\mathbb{G}^{(S,T]})$ predictable dual projection of $(\ind_{(S,T]})\stocint A$. Then, we have $(\ind_{(S,T]})\stocint A'' = (\ind_{(S,T]})\stocint A'$
\el

This lemma can be checked straightforwardly.

\bl\label{gSTKK}(Proposition \ref{appendix} (6))
For any $\mathbb{G}$-predictable process $K$, for any $(\mathbb{Q},\mathbb{G})$ local martingale $X$ such that $X_T\in\mathcal{G}^*_T$, the fact that $\ind_{(S,T]}K$ is $X$-integrable in $\mathbb{G}$ implies that $\ind_{(S,T]}K$ is $X^{(S,T]}$-integrable in $\mathbb{G}^{(S,T]}$. Conversely, for any $\mathbb{G}^{(S,T]}$-predictable process $K'$, for any $(\mathbb{Q},\mathbb{G}^{(S,T]})$ local martingale $X'$, the fact that $\ind_{(S,T]}K'$ is $X'$-integrable in $\mathbb{G}^{(S,T]}$ implies that $\ind_{(S,T]}K'$ is $X'^{(S,T]}$-integrable in $\mathbb{G}$. 
\el

\proof We note that the brackets $[X^{(S,T]}]$ and $[X'^{(S,T]}]$ are the same in the two filtrations. We can then check that $\sqrt{K\stocint[X^{(S,T]}]}$ or $\sqrt{K\stocint[X'^{(S,T]}]}$ is locally integrable in one filtration if and only if it is so in the other filtration.  \ok

We end the Appendix by the following result

\bl\label{GStar}(Lemma \ref{g-star})
Let $S,T$ be two $\mathbb{G}$-stopping times. We have$$
\mathcal{G}^*_{S\vee T}=\{S<T\}\cap \mathcal{G}^*_{T}+\{T\leq S\}\cap \mathcal{G}^*_{S}.
$$
\el

\proof Let us consider only this identity on the set $\{S<T\}\cap\{\tau\leq S\vee T<\infty\}$. With Lemma \ref{partTribu} we can write $$
\dcb
&&\{S<T\}\cap \{\tau\leq S\vee T<\infty\}\cap \mathcal{G}^*_{S\vee T}\\

&=&\{S<T\}\cap \{\tau\leq T<\infty\}\cap \sigma\{\{\tau\leq s\}, \{X_{S\vee T}\leq t\}: \mbox{ $X$ is $\mathbb{F}$ optional}, s,t\in\mathbb{R}\}\\

&=&\{S<T\}\cap \{\tau\leq T<\infty\}\cap \sigma\{\{\tau\leq s\}, \{X_{T}\leq t\}: \mbox{ $X$ is $\mathbb{F}$ optional}, s,t\in\mathbb{R}\}\\

&=&\{S<T\}\cap \{\tau\leq T<\infty\}\cap \mathcal{G}^*_T. \ok
\dce
$$
\end{document}